\documentclass{amsart}

\usepackage{egen,egenengelsk}
\DeclareMathOperator{\cok}{cok}
\DeclareMathOperator{\image}{im}
\DeclareMathOperator{\Exp}{Exp}
\DeclareMathOperator{\shExp}{Exp_{sh}}
\DeclareMathOperator{\locExp}{Exp_{loc}}

\DeclareMathOperator{\shLog}{Log_{sh}}
\DeclareMathOperator{\tr}{tr}
\newcommand{\loc}{\text{loc}}
\newcommand{\vir}{\text{vir}}

\newcommand{\Syst}{\mathbf{Syst}}
\newcommand{\SSyst}{\mathbf{SSyst}}

\usepackage[style=alphabetic, backend=biber, giveninits=true]{biblatex}
\DeclareFieldFormat{postnote}{#1} 

\author{Jørgen Vold Rennemo}
\address{Department of Mathematics\\
	University of Oslo\\
	P.O. Box 1053 Blindern\\
	0316 Oslo\\
	Norway}
\email{jorgeren@uio.no}


\title{Factorisability in K-theory}

\begin{document}
	\begin{abstract}
	We revisit and give a detailed proof of a lemma of Okounkov showing that, for a scheme $X$ with a torus action, the Euler characteristic generating function associated with a ``factorisable'' sequence of torus-equivariant coherent sheaves on the symmetric powers $\Sym^n X$ equals the plethystic exponential of the generating function of Euler characteristics of some sequence of sheaves on $X$.
	\end{abstract}
	
	\maketitle
	
	\setcounter{section}{1}
	
	Let $X$ be a complex, quasi-projective scheme equipped with an action of a torus $\mathbb{T} = (\CC^*)^d$, and let $(\cF_n)_{n=1}^\infty$ be a sequence where each $\cF_n$ is a coherent, $\mathbb{T}$-equivariant, cohomologically $\ZZ/2$-graded sheaf on $\Sym^n X$.
	We assume that $X^{\mathbb{T}}$ is proper, and so by localisation we may define $\chi(\Sym^n X, \cF_n)$ as a rational function in the coordinates of $\mathbb T$ (\cite[Thm.~2.2.]{thomason_lefschetz_1992}, see also \cite[2.3.10]{okounkov_lectures_2015}).
	
	Consider the addition morphism
	\begin{align*}
	\sigma_{m,n} \colon \Sym^{m}(X) \times \Sym^n(X) &\to \Sym^{m + n}(X)\\
	\left(\sum_{i=1}^m x_i, \sum_{i=1}^n y_i\right) &\mapsto \left(\sum_{i=1}^m x_i + \sum_{i=1}^n y_i\right).
	\end{align*}
	Let $U_{m,n} \subseteq \Sym^m(X) \times \Sym^n(X)$ be the open subset of pairs $(\sum_{i=1}^m x_i, \sum_{i=1}^n y_i)$ such that $x_i \not= y_j$ for all $i,j$.
	Following \cite[Sec.~5.3]{okounkov_lectures_2015}, we consider the notion of the sequence $(\cF_n)_{n=1}^\infty$ being \textit{factorisable}, which in particular means that 
	\begin{equation}
		\label{eqn:factorisationIsomorphism}
	\sigma_{m,n}^*(\cF_{m+n})|_{U_{m,n}} \cong (\cF_m \boxtimes \cF_n)_{U_{m,n}}
	\end{equation}
	as $\mathbb T$-equivariant sheaves for all $m,n$.
	Given positive integers $n_1,\dots, n_k$ with $n = \sum n_i$ and distinct points $x_1,\dots, x_k \in X$, then if $(\cF_n)_{n=1}^\infty$ is factorisable, at a point $y = \sum_{i=1}^k n_i x_i \in \Sym^{n}(X)$ we have
	\[
	\cF_n|_y \cong \bigotimes_{i=1}^k \cF_{n_i}|_{n_ix_i}.
	\]
	The full definition of factorisability (Def.\ \ref{def:factorisableSequence}) requires moreover that the isomorphisms \eqref{eqn:factorisationIsomorphism} are compatible in two ways representing associativity and commutativity.
	
	Easy examples of factorisable sequences include $\cF_n = \cO_{\Sym^n(X)}$, and 
	\[
	\cF_n = R\pi_*(\oplus_{i=0}^n \wedge^i \cO_X^{[n]}[i]),
	\]
	where $\cO_X^{[n]}$ the tautological rank $n$ bundle on $\Hilb^n(X)$, and $\pi \colon \Hilb^n X \to \Sym^n X$ is the Hilbert--Chow morphism.
	Harder examples of factorisable sequences include the direct images of virtual and twisted virtual structure sheaves on $\Hilb^n(X)$ for a Calabi--Yau 3-fold $X$, i.e.~$R\pi_*\cO_{\Hilb^n(X)}^{\vir}$ and $R\pi_*\widehat{\cO}_{\Hilb^n(X)}^{\vir}$ \cite{okounkov_lectures_2015, thimm_ktheoretic_2024}, as well as the sheaves on $\Sym^n(\CC^4)$ underlying the ``Nekrasov genus'' \cite{kool_proof_2025}.

Given a series $f(t_1,\dots, t_d,q) = \sum_{i=1}^\infty f_i(t_1,\dots,t_d)q^i$, where the $f_i$ are formal Laurent series or rational functions in the $t_i$, define the plethystic exponential by
\begin{equation}
\label{eqn:plethysticExp}
\Exp\left(f(t_1,\dots,t_d,q)\right) = \exp\left(\sum_{n \ge 1} \frac{f(t_1^n,\dots, t_d^n,q^n)}{n}\right).
\end{equation}
The goal of this paper is to show:
\begin{nprop}
	\label{thm:mainProposition}
	Let $X$ be a complex, quasi-projective, $\mathbb T$-equivariant scheme such that $X^\mathbb T$ is proper, and let $(\cF_n)_{n=1}^\infty$ be a sequence where each $\cF_n$ is a $\mathbb T$-equivariant, $\ZZ/2$-graded coherent sheaf on $\Sym^n X$.
	If $(\cF_n)_{n=1}^\infty$ is factorisable in the sense of Definition \ref{def:factorisableSequence}, then there exists a sequence $(\cE_n)^\infty_{n=1}$ of $\mathbb T$-equivariant, $\ZZ/2$-graded coherent sheaves on $X$ such that
	\[
	1 + \sum_{n=1}^\infty \chi(\Sym^n X,\cF_n)q^n = \Exp\left(\sum_{n=1}^\infty \chi(X, \cE_n)q^n\right).
	\]
\end{nprop}
This result (in a geometric form essentially equivalent to Proposition \ref{thm:mainGeometricProposition}) was stated by Okounkov as \cite[Lemma 5.3.4]{okounkov_lectures_2015}, though we found it hard to fill in the details in the proof given.
Our contribution here is to supply a complete and detailed proof of the statement, mainly with a view to applications in \cite{kool_proof_2025}, and to a generalised version of the proposition which is shown and used in \cite{thimm_ktheoretic_2024}.

Note that it is always true that a series $f(t_1,\dots, t_d,q)$ as on the left hand side of the proposition can be written as a plethystic exponential of some other formal series.
From an enumerative geometry perspective, the importance of Proposition \ref{thm:mainProposition} is that the coefficients $\chi(X,\cE_n)$ of the series on the right hand side typically have a special form as a rational function in the $t_i$.\footnote{E.g.~if $X = \AA^d$ with the standard $\mathbb T$-action, then $\chi(\AA^d,\cE_n)$ is of the form $g\prod_{i=1}^d(1-t_i)^{-1}$, where $g$ is a Laurent polynomial in the $t_i$.}
Thus the proposition provides strong constraints on the series on the left hand side; this plays a key role in various computations of K-theoretic DT3 and DT4 invariants \cite{kool_proof_2025,thimm_ktheoretic_2024,okounkov_lectures_2015}.

\subsection{Sketch of proof}
We now outline the proof of Proposition \ref{thm:mainProposition}, leaving all precise definitions and results to later sections.

Given a sequence of sheaves $(\cF_n)_{n=1}^\infty$, each $\cF_n$ living on $\Sym^n X$, we obtain a sequence of sheaves $F_n = p^*\cF_n$ on $[X^n/S_n]$ by the morphism $p \colon [X^n/S_n] \to \Sym^n X$.
From the $F_n$, we can compute $\chi(\Sym^n X,\cF_n) = \chi(X^n,F_n)^{S_n}$.
Working $S_n$-equivariantly on $X^n$ proves more convenient for our purposes than working with $\Sym^n X$, and in most of the paper we will work with a sequence $(F_n)_{n=1}^\infty$ where $F_n$ is a sheaf on $[X^n/S_n]$.

If $f(t_1,\dots, t_d, q)$ is the character of a $\ZZ/2$-graded $\mathbb T \times \CC^*_q$-representation $V$, then the plethystic exponential $\Exp(f)$ is the character of the infinite symmetric product of $V$.
In Section \ref{sec:plethysticBackground}, we study the analogous infinite symmetric product operation on sheaves on $\sqcup_{n=1}^\infty [X^n/S_n]$, which we call the sheaf-based plethystic exponential and denote by $\shExp$.
For a set $A \subseteq \{1,\dots, n\}$, let $i_A \colon X^n \to X^{|A|}$ be the natural projection $(x_1,\dots, x_n) \mapsto (x_i)_{i \in A}$, and let $F_A = i_A^*F_{|A|}$.
For a set partition $\mathfrak A = \{A_1,\dots, A_k\}$ of $\{1,\dots, n\}$, we let $F_{\mathfrak A} = \bigotimes_{i=1}^k F_{A_i}$.
We then define
\[
\shExp\left(\sum_{n=1}^\infty F_nq^n\right) = 1 + \sum_{n=1}^\infty \left(\bigoplus_{\mathfrak A} F_{\mathfrak A}\right) q^n,
\]
where the direct sum is over all set partitions $\mathfrak A$ of $\{1,\dots, n\}$, and this direct sum is equipped with a natural $S_n$-equivariant structure to make it a sheaf on $[X^n/S_n]$.

The operation $\shExp$ has a K-theoretic inverse denoted $\shLog$, defined order by order in $n$ so that
\[
\shLog\left(\shExp\left(\sum_{n=1}^\infty F_nq^n\right)\right) = \sum_{n=1}^\infty (F_n \oplus J_n \oplus J_n[1]) q^n
\]
for some sheaves $J_n$.
The operation $\shLog$ is not an inverse to $\shExp$ on the level of sequences of sheaves, but it is so at the level of K-theory.
The claim of Proposition \ref{thm:mainProposition} can then be reduced (see Section \ref{sec:proof}) to Proposition \ref{thm:mainGeometricProposition}, which says that if $(F_n)_{n=1}^\infty$ arises from a factorisable sequence on $\Sym^n X$, then
\[
\shLog\left(\sum_{n=1}^\infty F_nq^n\right) = \sum_{n=1}^\infty G_nq^n
\]
where each $G_n$ is equal in K-theory to a sheaf supported on the small diagonal $\Delta \subseteq X^n$; equivalently, that $[G_n]|_{X^n \setminus \Delta} = 0$.

In Section \ref{sec:explicitPlethysticLogarithm}, we pass from the abstract order-by-order definition of $\shLog$ to an explicit description of the coefficients $G_n$ as an sum of shifted copies of $F_{\mathfrak A}$, indexed by ``index trees of order $n$''.

If the $F_n$ come from a factorisable sequence of sheaves on $\Sym^n X$, then the sheaves $F_{\mathfrak A}$ on $X^n$ come equipped with partially defined isomorphisms between them: Given set partitions $\mathfrak A$ and $\mathfrak B$ of $\{1,\dots, n\}$ such that $\mathfrak A$ is finer than $\mathfrak B$, then $F_{\mathfrak A}$ is naturally isomorphic to $F_{\mathfrak B}$ on the set
\[
U_{\mathfrak A, \mathfrak B} = \{(x_1,\dots, x_n) \in X^n \mid \text{ if }i \not\sim_{\mathfrak A} j \text{ and } i \sim_{\mathfrak B} j, \text{ then }x_i \not= x_j\}.
\]
Furthermore these isomorphisms are compatible under compositions and with the natural $S_n$-actions.
In Section \ref{sec:systems}, we abstract the properties of such collections of sheaves and partial isomorphisms into the notion of a ``system of sheaves''.
The main result of this section is Corollary \ref{thm:KTheoryMapSurjective}, which says that every system of sheaves is equal in K-theory to a ``strict'' system, meaning one where the partial isomorphisms can be extended compatibly to global homomorphisms.

Applying this strictification, we replace the $F_{\mathfrak A}$ with K-theoretically equivalent sheaves, where the strictness property gives that for all $\mathfrak A, \mathfrak B$ with $\mathfrak A$ finer than $\mathfrak B$ we have homomorphisms $\phi_{\mathfrak A,\mathfrak B} \colon F_\mathfrak A \to F_\mathfrak B$ restricting to isomorphisms on $U_{\mathfrak A,\mathfrak B}$.
After making the corresponding replacement to $G_n$, in Section \ref{sec:definingAComplex} we then use these homomorphisms to put a differential $d$ on $G_n$.
We furthermore show that the cohomology sheaf $H^*(G_n, d)$ vanishes on $X^n \setminus \Delta$, and so $[G_n]|_{X^n \setminus \Delta} = 0$.
It follows that $[G_n] = [H_n]$, for some sheaf $H_n$ supported on $\Delta$, which is what we wanted.

\subsubsection{The proof for $n=2, 3$}
Let $(\cF_n)_{n=1}^\infty$ be a factorisable sequence of sheaves on $\Sym^n X$, which pulls back to a sequence $(F_n)_{n=1}^\infty$ on $[X^n/S_n]$.
Let $G_n$ be the coefficient of $q^n$ in $\shLog(\sum F_n q^n)$.
The central claim to show is Proposition \ref{thm:mainGeometricProposition}, that the sheaf $G_n$ is equal in K-theory to a sheaf $H_n$ pushed forth from the small diagonal $\Delta \subseteq X^n$.
As the general proof is technical, in this section we sketch how the proof goes for $n = 2$ and $3$, hoping to showcase the main ideas.

For $n = 2$, we have $G_2 = F_2 \oplus F_{\{1\},\{2\}}[1]$, where $F_{\{1\},\{2\}} = F_1 \boxtimes F_1$.
Let $U = X^2 \setminus \Delta$, and let $j \colon U \to X^2$ be the inclusion.
The factorisability condition on the $F_i$ implies that there exists an $S_2$-invariant isomorphism
\[
\phi \colon F_{\{1\},\{2\}}|_U \to F_2|_U.
\]
Thus the class $[G_2] = [F_2] - [F_{\{1\},\{2\}}]$ vanishes on $U$, and hence equals the class of some sheaf $H_2$ supported on $\Delta$ by the localisation sequence in K-theory of coherent sheaves (see e.g.~\cite[2.2.19]{okounkov_lectures_2015} for an explicit argument).

In the more complicated case of $n > 2$, where $G_n$ is a sum of multiple $F_{\mathfrak A}$ which are identified on different open subsets of $X^n$, we are not able to use the localisation sequence in this way.
Our general proof of the vanishing of $[G_n]|_U$ specialises inthe case of $n = 2$ as follows:
First, the strictification procedure of Section \ref{sec:systems} replaces $F_2$ and $F_{\{1\},\{2\}}$ by K-theoretically equivalent sheaves.
We have homomorphisms
\[
a \colon F_{\{1\},\{2\}} \to j_*j^*F_{\{1\},\{2\}}
\]
\[
b \colon F_{\{1\},\{2\}} \to j_*j^*F_{\{1\},\{2\}} \overset{j_*\phi}{\to} j_*j^* F_2 
\]
\[
c \colon F_2 \to j_*j^*F_2.
\]
Setting
\[
F_{\{1\},\{2\}}' = \im a \oplus \ker a,
\]
and
\[
F_2' = (\im b + \im c) \oplus \ker c \oplus (\im b + \im c)/\im c[1],
\]
we have $[F_{\{1\},\{2\}}] = [F'_{\{1\},\{2\}}]$ and $[F_2] = [F'_2]$.
Moreover, we have a homomorphism
\[
\phi' \colon F'_{\{1\},\{2\}} \to \im a \overset{j_*j^*\phi}{\to} \im b \to \im b + \im c \to F'_2.
\]
It is easy to see that after restriction to $U$, we have $F'_{\{1\},\{2\}} = F_{\{1\},\{2\}}, F'_2 = F_2$, and $\phi' = \phi$, in particular $\phi'$ is an isomorphism over $U$.
The complex $G_2' = F'_{\{1\},\{2\}} \overset{\phi'}{\to} F'_2$ then has cohomology supported (set-theoretically) on $\Delta$, and so by devissage there exists an $H_2$ supported on $\Delta$ such that $[G_2] = [G'_2] = [H_2]$.

For the case of $n = 3$, we have (cf.~ Lemma \ref{thm:formulaForLog}) that
\[
G_3 = F_{3} \oplus (F_{\{\{2,3\},\{1\}\}} \oplus F_{\{\{1,3\},\{2\}\}} \oplus F_{\{\{1,2\},\{3\}\}} \oplus F_{\{\{1\},\{2\},\{3\}\}})[1] \oplus F_{\{1\},\{2\},\{3\}}^{\oplus 3}[2]
\]
Here $F_{\{1\},\{2\},\{3\}} = F_1 \boxtimes F_1 \boxtimes F_1$, and $F_{\{i,j\},\{k\}} = \pi_{ij}^*F_2 \otimes \pi_k^*F_1$, with $\pi_{ij} \colon X^3 \to X^2$ and $\pi_k \colon X^3 \to X$ denoting the projections to factors $i,j$ and $k$, respectively.
The factorisability of the sequence $F_n$ gives us isomorphisms $\phi_{\mathfrak A,\mathfrak B}$ for any two set partitions $\mathfrak A, \mathfrak B$ of $\{1,2,3\}$ with $\mathfrak B$ refining $\mathfrak A$.
Concretely, this gives partially defined isomorphisms
\[
\phi_k \colon F_{\{1\},\{2\},\{3\}} \to F_{\{i,j\},\{k\}} \text{ on the locus where }x_i \not= x_j,
\] 
\[
\psi_k \colon F_{\{i,j\},\{k\}} \to F_3 \text{ on the locus where }x_i, x_j \not= x_k,
\]
and
\[
\omega \colon F_{\{1\},\{2\},\{3\}} \to F_3 \text{ on the locus where all the }x_i\text{ are distinct.}
\]
Moreover 
\begin{equation}
	\label{eqn:associativity}
	\omega = \psi_1 \circ \phi_1 = \psi_2 \circ \phi_2 = \psi_3 \circ \phi_3
\end{equation}
on the locus where all $x_i$ are distinct.

The strictification procedure of Section \ref{sec:systems} allows us to replace the $F_{\mathfrak A}$ by sheaves with the same K-theory class, with the property that the partial isomorphisms $\phi_k,\psi_k,\omega$ extend to homomorphisms on the whole of $X^3$ in such a way that \eqref{eqn:associativity} holds globally.
After this replacement, we may endow $G_3$ with a differential defining the complex
\[
G_3 = F_{\{1\},\{2\},\{3\}}^{\oplus 3} \to F_{\{\{2,3\},\{1\}\}} \oplus F_{\{\{1,3\},\{2\}\}} \oplus F_{\{\{1,2\},\{3\}\}} \oplus F_{\{\{1\},\{2\},\{3\}\}} \to F_{3},
\]
Here the first differential is 
\[
\begin{pmatrix} -\phi_1 & 0 & 0 \\ 
	0 & -\phi_2 & 0 \\
	0 & 0 & -\phi_3 \\
	\id & \id & \id 
\end{pmatrix}
\]
and the second is
\[
\begin{pmatrix} -\psi_1 & -\psi_2 & -\psi_3 & -\omega \end{pmatrix}.
\]
If $\bar x = (x_1,x_2,x_3) \in X^3 \setminus \Delta$, then without loss of generality we may assume that $x_1 \not= x_2, x_3$, and so $\psi_1, \phi_2$ and $\phi_3$ are all isomorphisms at $\bar x$.
The complex $G_3$ admits a filtration (cf.~Section \ref{sec:acyclicity}) with the following subquotients
\begin{gather*}
	F_{\{2,3\}, \{1\}} \overset{-\psi_1}{\to} F_3 \\
	F_{\{1\},\{2\},\{3\}} \overset{-\phi_2}{\to} F_{\{1,3\}, \{2\}} \\
	F_{\{1\},\{2\},\{3\}} \overset{-\phi_3}{\to} F_{\{1,2\}, \{3\}} \\
	F_{\{1\},\{2\},\{3\}} \overset{\id}{\to} F_{\{1\}, \{2\},\{3\}},
\end{gather*}
all of which are acyclic at $\bar x$.
As $\bar x$ was chosen arbitrarily, it follows that $G_3$ is acyclic on $X^3 \setminus \Delta$.
As the cohomology of $G_3$ is supported set-theoretically on $\Delta$, it follows by devissage that $[G_3] = [H_3]$ for some sheaf $H_3$ on $\Delta$.

\begin{nremark}
	The main difficulties in the proof arise from the fact that K-theory classes of coherent sheaves are not ``stratifiable'', in the sense that for a sheaf $F$ on a scheme $X$ with open subset $U$, there is no way to reconstruct $[F]$ from $[F|_{U}]$ and $[F|_{X \setminus U}]$.
	
	By contrast, if $F$ is a constructible sheaf, $j \colon U \to X$ and $i \colon X \setminus U \to X$ the inclusions, then we have $[F] = j_*[F|_U] + i_![F|_{X \setminus U}]$.
	The statement analogous to Proposition \ref{thm:mainGeometricProposition} for constructible sheaves can therefore be proved more simply by stratifying $X^n$ by the $S_n$-orbits of its diagonals and then checking that the restriction of $\shExp \left(\sum_{n=1}^\infty F_nq^n\right)$ to each stratum except for the small diagonal vanishes.
\end{nremark}

\subsection{Acknowledgements}
This paper was originally conceived as a paragraph -- then a few pages -- then an appendix -- in the author's joint paper \cite{kool_proof_2025} with Martijn Kool.
Conversations with Kool and with Felix Thimm have been essential throughout the writing of the paper.

The work has been funded by Research Council of Norway grant no.\ 302277.

The author acknowledges that Proposition \ref{thm:mainProposition} deserves a more natural proof, but what can you do.

\subsection{Conventions}
We let $[n] = \{1,\ldots, n\}$.

The power set $2^{[n]}$ is equipped with the ``binary'' total order such that for $A, B \in 2^{[n]}$, we have $A \le B$ if and only if
\[
\sum_{i \in A} 2^i \le \sum_{i \in B} 2^i,
\]
so in particular $A \subseteq B$ implies $A \le B$.

A set partition of a set $S$ is a set $\mathfrak A = \{A_1,\ldots, A_k\}$, where the $A_i$ are disjoint subsets of $S$ whose union is $S$.
Given a set partition $\mathfrak A = \{A_1, \ldots, A_k\}$ of $[n]$, we write $\underline{\mathfrak A}$ for the partition of $n$ given by the multiset $\{|A_i|\}_{i=1}^k$.

Given a series expression $f(q) = \sum_{n \in \ZZ}a_n q^n$, we let $\langle q^n\rangle (f(q))= a_n$.

All group actions are on the left.

\section{Plethystic exponentials and logarithms}
\label{sec:plethysticBackground}
	Let $X$ be a complex, quasi-projective scheme equipped with an action of a torus $\mathbb{T} = (\CC^*)^d$.
	We let $(F_n)_{n=1}^\infty$ be a sequence where $F_n$ is a coherent, $\mathbb{T} \times S_n$-equivariant, cohomologically $\ZZ/2$-graded sheaf on $X^n$.
	For the remainder of the paper, we will write ``sheaf'' for a coherent $\mathbb T$-equivariant $\ZZ/2$-graded sheaf, and so for instance say that $F_n$ is a sheaf on $[X^n/S_n]$ or a $S_n$-equivariant sheaf on $X^n$.

	For any finite set $A$, let $X^A$ be the direct product of copies of $X$, indexed by $A$. 
	By definition, $X^n = X^{[n]}$.
	The assignment $A \mapsto X^A$ defines a contravariant functor from finite sets to schemes.
	Given a map of sets $u \colon A \to B$, we will write $u_X \colon X^B \to X^A$ for the associated morphism of schemes.
	The group $S_n$ acts on $X^n$ by sending $\sigma \in S_n$ to $\sigma_X^{-1}$, so that
	\[
	\sigma(x_1,\dots,x_n) = (x_{\sigma^{-1}(1)},\dots, x_{\sigma^{-1}(n)}).
	\]
	
	Let $A$ be a subset of $[n]$.
	There is a unique order preserving inclusion $i_A \colon [|A|] \to [n]$ with image $A$, giving rise to a morphism $i_{A,X} \colon X^n \to X^{|A|}$, and we define
	\[
	F_{A} = i_{A,X}^*F_{|A|}.
	\]
	Let $\sigma \in S_n$.
	There is a unique $\tau \in S_{|A|}$ such that 
	\[
	\sigma^{-1} \circ i_A = i_{\sigma^{-1}(A)} \circ \tau^{-1},
	\]
	and so
	\[
	i_{A,X} \circ \sigma_X^{-1} = \tau_X^{-1} \circ i_{\sigma^{-1}(A),X}
	\]
	Since $F_{|A|}$ is $S_{|A|}$-equivariant, there is a given isomorphism
	\[
	\rho_{|A|, \tau} \colon (\tau^{-1}_X)^*F_{|A|} \cong F_{|A|},
	\]
	which gives rise to an isomorphism
\begin{align*}
	\label{eqn:basicIsomorphism}
	\rho_{A,\sigma} \colon (\sigma_X^{-1})^*F_A &= (\sigma_X^{-1})^*i_{A,X}^*F_{|A|} \\
	&=  i_{\sigma^{-1}(A),X}^*(\tau_X^{-1})^*F_{|A|} \overset{\rho_{|A|,\tau}}{\to} i_{\sigma^{-1}(A),X}^*F_{|A|} = F_{\sigma^{-1}(A)}.
\end{align*}
	
	Let
	\[
	\mathfrak A = \{A_1, \ldots, A_l\}
	\] 
	be a set partition of $[n]$, and assume that the $A_i$ are lexicographically ordered with $A_1 < \dots < A_l$.
	We define 
	\[
	F_{\mathfrak A} = F_{A_1} \otimes \cdots \otimes F_{A_l}.
	\]
	
	Using the isomorphisms $\rho_{A_i,\sigma}$ and the commutativity of the tensor product, we get an isomorphism
	\begin{align*}
    \rho_{\mathfrak A, \sigma} \colon (\sigma^{-1}_X)^*F_{\mathfrak A} &\cong (\sigma^{-1}_X)^*(F_{A_1}) \otimes \cdots \otimes (\sigma^{-1}_X)^*(F_{A_l}) \cong F_{\sigma^{-1}(A_1)} \otimes \cdots \otimes F_{\sigma^{-1}(A_l)} \\
    &\cong F_{\sigma^{-1}(A_{i_1})} \otimes \cdots \otimes F_{\sigma^{-1}(A_{i_l})} = F_{\sigma^{-1}(\mathfrak A)},
	\end{align*}
	where $i_1,\ldots, i_l$ are such that the $\sigma^{-1}(A_{i_j})$ forms an increasing sequence in the binary ordering of subsets of $[n]$.

	Let $\lambda$ be a partition of $n$. We define
	\[
	F_\lambda = \bigoplus_{\underline{\mathfrak A} = \lambda} F_{\mathfrak A}.
	\]
	This has the structure of an $S_n$-equivariant sheaf on $X^n$, where the $S_n$-structure map $\rho_{\lambda, \sigma} \colon (\sigma^{-1}_X)^* F_\lambda \to F_\lambda$ is the sum of the isomorphisms $\rho_{\mathfrak A, \sigma}$.
	
	\begin{nremark}
		\label{remark:alternativeDef}
		The sheaf $F_\lambda$ can alternatively be described as follows.
		Choose a set partition $\mathfrak A = \{A_1,\ldots, A_l\}$ of $[n]$ such that $\underline{\mathfrak{A}} = \lambda$.
		Let $S_{\mathfrak A} \subseteq S_n$ be the subgroup preserving $\mathfrak A$.
		The sheaf
		\[
		F^{\text{pre}}_\lambda = F_{A_1} \otimes \cdots \otimes F_{A_l}
		\]
		is then naturally equivariant with respect to the action of $S_{\mathfrak A}$ on $X^n$.
		There is a natural morphism of stacks $[X^n/S_{\mathfrak A}] \to [X^n/S_n]$, and $F_\lambda$ is isomorphic to the pushforward of $F^{\text{pre}}_\lambda$ under this morphism.
	\end{nremark}
	
	\begin{nremark}
		\label{remark:noncommutativity}
		In general, for sheaves $F$ and $G$, the diagram of isomorphisms
		\[
		\begin{tikzcd}
			F[i] \otimes G[j] \arrow[r] \arrow[d] & (F \otimes G)[i + j] \arrow[d] \\
			G[j] \otimes F[i] \arrow[r] & (G \otimes F)[i + j]
		\end{tikzcd}
		\]
		only commutes up to a Koszul sign $(-1)^{ij}$.
		Since the permutation of tensor factors is used in the definition of the $S_n$-equivariant structure of $F_\lambda$, replacing the sequence $(F_i)$ by some sequence $(F_i[k_i])$ will potentially change $F_\lambda$ more than by just a shift.
		
		As an example, take $\lambda = (1,1)$ and $F'_1 = F_1[1]$.
		We then have that $F'_{(1,1)} \cong F_{(1,1)} \otimes \mathrm{sgn}$, where $\mathrm{sgn}$ is the non-trivial character of $S_2$.
	\end{nremark}
	
	\subsection{The plethystic exponential on $X^n$}
	We define a plethystic exponential operation
	\[
	\shExp\left(\sum_{n=1}^\infty F_n q^n\right) = 1 + \sum_{n=1}^\infty \left(\bigoplus_{\lambda \vdash n} F_\lambda\right) q^n.
	\]
	We define a plethystic logarithm operation by
\begin{equation}
	\label{eqn:sheafLogarithm}
\shLog\left(1 + \sum_{n=1}^\infty F_nq^n\right) = \sum_{n=1}^\infty G_n q^n,
\end{equation}
where $G_n$ is defined inductively by
\[
G_n = F_n \oplus \langle q^n\rangle\left(\shExp\left(\sum_{k=1}^{n-1}G_kq^k\right)\right)[1],
\]

\subsection{Plethystic operations in K-theory}
The operations of plethystic exponential and logarithm defined above have input and output power series in $q$ whose $q^n$-coefficients are isomorphism classes of sheaves on $[X^n/S_n]$.

Given a sheaf $E$ on a stack $Y$, with $\ZZ/2$-graded pieces $E_{\text{even}}$ and $E_{\text{odd}}$, we define
\[
[E] = [E_{\text{even}}]-[E_{\text{odd}}] \in K(Y).
\]
For every partition $\lambda$ of $n$, the map $(F_n)_{n=1}^\infty \mapsto F_\lambda$ defined above factors through K-theory in the sense that $[F_\lambda] \in K_{\mathbb{T} \times S_n}(X^n)$ only depends on the classes $[F_i] \in K_{\mathbb{T} \times S_i}(X^i)$.

The plethystic exponentials and logarithms therefore descend to operations on power series with coefficients in $K_{\mathbb{T} \times S_n}(X^n)$.
Note that the plethystic exponentials and logarithms are not mutually inverse when thought of as operations on sheaves, but are so as operations on K-theory classes.

\subsection{Passing to localised K-theory}
Let $K_\mathbb{T}(\pt)_{\text{loc}}$ be the fraction field of the ring $K_{\mathbb{T}}(\pt)$.
If $\mathbb{T} = \Spec \CC[t^\pm_1,\ldots, t_d^{\pm}]$, then we have canonical isomorphisms $K_\mathbb{T}(\pt) \cong \ZZ[t_1^\pm,\cdots, t_d^{\pm}]$ and $K_\mathbb{T}(\pt)_\loc \cong \QQ(t_1,\ldots, t_d)$.
For a $K_\mathbb{T}(\pt)$-module $M$, let $M_{\text{loc}} = M \otimes_{K_\mathbb{T}(\pt)} K_\mathbb{T}(\pt)_{\text{loc}}$.

Given a sequence $(\alpha_i)_{i=1}^\infty$, where $\alpha_i \in K_{\mathbb{T} \times S_i}(X^i)_{\text{loc}}$, we want to define a plethystic exponential
\begin{equation}
	\label{eqn:plethysticExpKTheory}
\locExp\left(\sum_{i=1}^\infty \alpha_iq^i \right) = 1 + \sum_{i=1}^\infty \beta_iq^i
\end{equation}
with $\beta_i \in K_{\mathbb{T} \times S_i}(X^i)_{\text{loc}}$, agreeing with $\shExp$ for series with coefficients in $K_{\mathbb{T} \times S_i}(X^i)$.

We first extend the definition of the operations $F_\lambda$, where $\lambda$ is a partition of $n$.
\begin{nlemma}
	\label{thm:ElambdaInvertible}
	Let $D$ be a sheaf\footnote{Recall this means a $\mathbb{T}$-equivariant $\ZZ/2$-graded coherent sheaf on $\pt$.} on $\pt$.
	For each $i \ge 1$, let $D_i = D$, considered as a sheaf on $[\pt/S_i]$ with the trivial $S_i$-equivariant structure.
	
	Let $\mathfrak A = \{A_1,\ldots, A_l\}$ be a set partition of $[n]$, and let $S_{\mathfrak A} \subseteq S_n$ be the subgroup fixing $\mathfrak A$.
	If $[D] \not= 0 \in K_{\mathbb{T}}(\pt)$, then
	\[
		[D_{A_1} \otimes \cdots \otimes D_{A_l}]
	\]
	is a unit in $K_{\mathbb{T} \times S_{\mathfrak A}}(\pt)_{\mathrm{loc}}$.
\end{nlemma}
\begin{proof}
	Let $\epsilon = [D_{A_1} \otimes \cdots \otimes D_{A_l}]$.
	The $K_\mathbb{T}(\pt)_{\text{loc}}$-vector space $K_{\mathbb{T} \times S_{\mathfrak A}}(\pt)_{\text{loc}}$ is finite-dimensional, so it is enough to see that $\epsilon$ is not a 0-divisor (an inverse is then given by some polynomial in $\epsilon$ by the Cayley--Hamilton theorem).
	
	Let $\sigma \in S_{\mathfrak A}$, let $\overline{\sigma}$ be the image of $\sigma$ under the natural homomorphism $S_{\mathfrak A} \to S_l$, and let $(c_1,\cdots,c_r)$ be the lengths of the cycles in $\overline{\sigma}$.
	For any $t \in \mathbb{T}$, one computes
	\[
	\tr_{(t,\sigma)}(\epsilon) = \prod_{i=1}^r\tr_{t^{c_i}}(D).
	\]
	In particular, since $[D] \not= 0$, we have that $\tr_{t}(D)$ is not 0 as a function of $t$, so $\tr_{(t,\sigma)}(\epsilon)$ is not 0 as a function of $t$.
	
	Now if $0 \not= \alpha \in K_{\mathbb{T} \times S_{\mathfrak A}}(\pt)_{\text{loc}}$, there is some $(t,\sigma) \in \mathbb{T} \times S_{\mathfrak A}$ such that $\tr_{(t,\sigma)}(\alpha) \not= 0$, and so
	\[
	\tr_{(t,\sigma)}(\epsilon\alpha) = \tr_{(t,\sigma)}(\epsilon)\tr_{(t,\sigma)}(\alpha) \not= 0,
	\]
	hence $\epsilon\alpha \not= 0$.
\end{proof}

Given classes $\alpha_1,\ldots, \alpha_n$ with $\alpha_i \in K_{\mathbb{T} \times S_i}(X^i)_{\text{loc}}$, we can find sheaves $F_i$ on $[X^i/S_i]$ and $D$ on $\pt$ such that 
\[
\alpha_i = \frac{[F_i]}{[D]}
\]
for all $i$.
For a partition $\lambda$ of $n$, we then define $\alpha_\lambda$ as follows.
Let $\mathfrak A = \{A_1,\ldots, A_l\}$ be a set partition of type $\lambda$, and let $p \colon [X^n/S_{\mathfrak A}] \to [X^n/S_n]$ be the natural morphism.
We let
\[
\alpha_\lambda = p_*([D_{A_1}\otimes \cdots \otimes D_{A_l}]^{-1}[F_{A_1} \otimes \cdots \otimes F_{A_l}]).
\]
By Lemma \ref{thm:ElambdaInvertible}, the inverse in the expression above exists.
One checks that it does not depend on the chosen presentation of $\alpha_i$ as a quotient, using the fact that given sheaves $(F_1,\ldots, F_n)$ and $(G_1,\ldots, G_n)$, with $H_i = F_i \otimes G_i$, we have
\[
H_{A_1} \otimes \cdots \otimes H_{A_l} \cong (F_{A_1}\otimes \cdots \otimes F_{A_l}) \otimes (G_{A_1} \otimes \cdots \otimes G_{A_l})
\]
as sheaves on $[X^n/S_{\mathfrak A}]$.
Finally, by Remark \ref{remark:alternativeDef}, taking $[D] = 1$ shows that if $\alpha_i = [F_i]$, then $\alpha_\lambda = [F_\lambda]$.

The plethystic exponential in localised $K$-theory is then defined by
\[
\locExp\left(\sum_{n=1}^\infty \alpha_n q^n\right) = 1 + \sum_{n=1}^\infty \sum_{\lambda \vdash n} \alpha_\lambda q^n
\]

The following lemma shows that the operation of $\locExp$ agrees with $\Exp$ as in \eqref{eqn:plethysticExp} under taking characters.
\begin{nlemma}
	\label{thm:numericalLemma}
	Let $(\alpha_i)_{i=1}^\infty$ be a sequence with $\alpha_i \in K_{\mathbb{T} \times S_i}(\pt)_{\loc}$, and let $\alpha_i^{S_i} \in K_{\mathbb{T}}(\pt)_{\loc}$ be the invariant part.
	Let $1 + \sum_{i \ge 1} \beta_iq^i = \locExp\left(\sum_{i \ge 1} \alpha_iq^i\right)$.
	Then we have an equality in $\QQ(t_1,\ldots, t_d)[[q]]$
	\begin{equation}
		\label{eqn:lemmaEquation}
		1 + \sum_{n=1}^\infty \tr_t(\beta_n^{S_n})q^n = \exp\left(\sum_{i=1}^\infty \frac{\sum_{n\ge 1}\tr_{t^i}(\alpha_n^{S_n})q^{in}}{i}\right)
	\end{equation}
\end{nlemma}
Note that by definition, the right hand side of \eqref{eqn:lemmaEquation} equals
\[
\Exp\left(\sum_{n=1}^\infty \tr_t(\alpha_n^{S_n})q^n\right)
\]
\begin{proof}
	Let $\lambda = (1^{l_1}\cdots k^{l_k})$ be a partition of $n = \sum_{i=1}^k il_i$.
	It is easy to check that 
	\[
	\alpha_\lambda^{S_n} = \alpha_{(1^{l_1})}^{S_{l_1}}\cdots \alpha_{(k^{l_k})}^{S_{kl_k}}.
	\]
	From this we get
	\begin{align*}
		1 + \sum_{n \ge 1} \beta_n^{S_n}q^n = 1 + \sum_{n \ge 1}\sum_{\lambda \vdash n} \alpha^{S_n}_\lambda q^n = \prod_{i \ge 1}\left(1 + \sum_{n \ge 1} \alpha_{(i^n)}^{S_{in}}q^{in}\right) = \prod_{i \ge 1} \locExp\left(\alpha_i q^i\right)^{S_*},
	\end{align*}
	where the superscript $S_*$ denotes replacing the $q^k$-coefficient by its $S_k$-invariants for all $k$.
	Since the right hand side of \eqref{eqn:lemmaEquation} is multiplicative under addition of $\alpha_i$, we reduce to the case where $\alpha_i \not= 0$ for a single $i$.
	
	We may write $\alpha_i = [F]/[D]$, where $F$ and $D$ are sheaves on $[\pt/S_i]$ and $D$ equals the pullback of some sheaf $D'$ on $\pt$.
	The only partitions $\lambda$ such that $\alpha_\lambda \not= 0$ are of the form $\lambda = (i^n)$, and so we have $\beta_{in} = \alpha_{(i^n)}$, with $\beta_k = 0$ for all $k$ not divisible by $i$.
	
	Let $\mathfrak A$ be a set partition of $[in]$ of type $(i^n)$.
	There is a natural commutative diagram
	\[
	\begin{tikzcd}
		{[\pt/S_{\mathfrak A}]} \arrow[r, "f"] \arrow[d, "p"] & {[\pt/S_{in}]} \arrow[d, "q"] \\
		{[\pt/S_{n}]} \arrow[r, "g"] & \pt,
	\end{tikzcd}
	\]
	and we have by definition
	\[
	\alpha_{(i^n)} = f_*([F_{\mathfrak A}][D_{\mathfrak A}]^{-1}).
	\]
	Taking $S_{in}$-invariants, we have
	\[
	\alpha_{(i^n)}^{S_{in}} = q_*f_*([F_{\mathfrak A}][D_{\mathfrak A}]^{-1}) = g_*p_*([F_{\mathfrak A}][D_{\mathfrak A}]^{-1}).
	\]
	Let now 
	\[
	\gamma = p_*([F_{\mathfrak A}][D_{\mathfrak A}]^{-1}).
	\]
	Since $D$ is the pullback of the sheaf $D'$ on $\pt$, we have
	\[
	D_{\mathfrak A} = p^*((D')^{\otimes n}),
	\]
	Noting $p_*(F_{\mathfrak A}) = (F^{S_i})^{\otimes n}$ and applying the projection formula, we get
	\[
	\gamma = p_*([F_{\mathfrak A}][D_{\mathfrak A}]^{-1}) = [(F^{S_i})^{\otimes n}][(D')^{\otimes n}]^{-1}
	\]
	
	For $t \in \mathbb{T}$, we have
	\[
	\tr_t(\beta_{in}^{S_{in}}) = \tr_t(g_*\gamma) = \frac{\sum_{\sigma \in S_n} \tr_{(t,\sigma)}(\gamma)}{n!}
	\]
	If $\sigma$ is of cycle type $(1^{a_1}\cdots k^{a_k})$, we have
	\[
	\tr_{(t,\sigma)}(\gamma) = \frac{\tr_{(t,\sigma)}((F^{S_i})^{\otimes n})}{\tr_{(t,\sigma)}((D')^{\otimes n})} = \frac{\prod_{j=1}^k\tr_{t^j}(F^{S_i})^{a_j}}{\prod_{j=1}^k \tr_{t^j}(D')^{a_j}}.
	\]
	Since 
	\[
	\alpha_i^{S_i} = [F^{S_i}]/[D'],
	\]
	we get 
	\[
	\tr_{(t,\sigma)}(\gamma) = \prod_{j=1}^k \tr_{t^j}(\alpha_i^{S_i})^{a_j}.
	\]
	There are
	\[
	\frac{n!}{\prod_{j=1}^k j^{a_j}a_j!}
	\]
	elements in $S_n$ of cycle type $(1^{a_1}\cdots k^{a_k})$, so we get
	\[
	\tr_t(\beta_{in}^{S_{in}}) = \sum_{\lambda = (i^{a_i}) \vdash n} \frac{\prod_{j=1}^k \tr_{t^j}(\alpha_i^{S_i})^{a_j}}{\prod_{j=1}^k j^{a_j}a_j!},
	\]
	which gives a computation of the left hand side of \eqref{eqn:lemmaEquation}.
	This can be matched to the result of expanding the exponential on the right hand side of \eqref{eqn:lemmaEquation}.
\end{proof}

\section{The plethystic logarithm formula}
\label{sec:explicitPlethysticLogarithm}
In this section, we give an explicit presentation of the operation $\shLog$ from \eqref{eqn:sheafLogarithm}, writing the components $G_n$ as a direct sum of shifted copies of $F_{\mathfrak A}$, equipped with a modified $S_n$-equivariant structure.
\begin{ndefn}
	An \textbf{index tree of order $n$} is the data of
	\begin{itemize}
		\item a rooted tree such that each node is either a leaf or has at least two children.
		\item for each leaf $v$ of the tree, a subset $L_v \subseteq [n]$, chosen such that
		\[
		\{L_v \mid v\text{ a leaf of }T\}
		\]
		is a partition of $[n]$.
	\end{itemize}
\end{ndefn}
We let $\mathcal T_n$ be the set of index trees of order $n$.
The nodes of an index tree $T$ are partially ordered by the relation $\preceq$, such that $v \preceq w$ if $v$ lies on the unique path connecting the root node and $w$.
We then say that $w$ is a descendant of $v$ and $v$ an ancestor of $w$.
For any non-leaf node $v$ of $T$, we let 
\[
L_v = \bigcup_{\substack{v' \succeq v \\
		v' \text{ a leaf}}} L_{v'}.
\]
\begin{nremark}
	\label{rmk:treesAsSets}
	Since $v \preceq v'$ if and only if $L_{v'} \subseteq L_v$, and since the tree structure in turn is determined by $\preceq$, an index tree $T$ is fully determined by the set
	\[
	L(T) = \{L_v \mid v\text{ a node of }T\} \subseteq 2^{[n]}.
	\]
	Using the correspondence $T \leftrightarrow L(T)$, an index tree is equivalent to the data of a subset $S \subseteq 2^{[n]}$ such that
	\begin{itemize}
		\item $[n] \in S$
		\item Given two elements $L_1, L_2 \in S$, either $L_1 \cap L_2 = \varnothing$ or one is contained in the other.
		\item Every element of $S$ is either minimal under inclusion or the union of its strict subsets in $S$.
	\end{itemize}
\end{nremark}

\begin{ndefn}
	For an index tree $T$ of order $n$, we let $\mathfrak A(T)$ be the partition of $[n]$ given by
	\[
	\mathfrak A(T) = \{L_v \mid v\text{ a leaf of }T\}
	\]
	We let 
	\[
	P(T) = L(T) \setminus \mathfrak A(T) = \{L_v \mid v\text{ a non-leaf node of }T\}
	\]
\end{ndefn}

\begin{ndefn}
	Define the sheaf $F_T$ on $X^n$ by
	\[
	F_T = F_{\mathfrak A(T)}.
	\]
\end{ndefn}

\begin{ndefn}
	Given an element $\sigma \in S_n$, and an index tree $T$, define $\sigma(T)$ to be the same underlying tree where the node $v$ is assigned the label set $\sigma(L_v)$.
\end{ndefn}

Let $T$ be an index tree of order $n$, and let $\sigma \in S_n$.
Working with the binary ordering on $P(T)$, let $l(T,\sigma)$ be the number of pairs $A < B \in P(T)$ such that $\sigma(A) > \sigma(B)$.
We define an isomorphism
\begin{equation}
	\label{eqn:rhoSigma}
\rho_{T, \sigma} \colon (\sigma^{-1}_X)^* F_T \to F_{\sigma^{-1}(T)}
\end{equation}
as $(-1)^{l(T,\sigma)}$ times the composition
\[
(\sigma^{-1}_X)^* F_T = (\sigma^{-1}_X)^*F_{\mathfrak A(T)} \overset{\rho_{\mathfrak A(T),\sigma}}\to F_{\sigma^{-1}\mathfrak A(T)} = F_{\mathfrak A(\sigma^{-1}(T))} = F_{\sigma^{-1}(T)}.
\]

Let $T$ be an index tree of order $n$, with $S_n$-orbit $S_nT$ under the action of $S_n$ on the set of index trees $\mathcal T_n$.
Define a sheaf $F_{S_nT}$ on $[X^n/S_n]$ by
\[
F_{S_nT} = \bigoplus_{T' \in S_nT} F_{T'},
\]
where we give $F_{S_nT}$ the $S_n$-equivariant structure defined by the sum of the $\rho_{T,\sigma}$ from \eqref{eqn:rhoSigma}.

\begin{nremark}
	\label{remark:alternativeDefT}
	Similarly to Remark \ref{remark:alternativeDef}, we could equivalently define $F_{S_nT}$ as follows.
	Let $S_T \subseteq S_n$ be the subgroup fixing $T$.
	Then $S_T \subseteq S_{\mathfrak A(T)}$, and so we have a diagram
	\[
	\begin{tikzcd}
		{[X^n/S_T]} \arrow[r, "q"] \arrow[d, "p"] & {[X^n/S_n]} \\
		{[X^n/S_{\mathfrak A(T)}]} &
	\end{tikzcd}
	\]
	The action of $S_T$ on $P(T)$ gives a sign character $\mathrm{sgn}_{P(T)}$ of $S_T$, and we have 
	\[
	F_{S_nT} \cong q_*(p^*(F_\mathfrak A(T)) \otimes \mathrm{sgn}_{P(T)}).
	\]
\end{nremark}

For every $k \ge 0$, we let $\cT_{n,k} \subseteq \cT_n$ be the set of index trees of order $n$ with $k$ non-leaf nodes.
This is non-empty precisely when $0 \le k \le n-1$.
Let
\begin{equation}
	\label{eqn:definitionOfGnk}
G_{n,k} = \bigoplus_{S_nT \in \mathcal T_{n,k}/S_n} F_{S_nT},
\end{equation}
\begin{nlemma}
	\label{thm:formulaForLog}
	Let $(F_n)_{n=1}^\infty$ be a sequence where $F_n$ is a sheaf on $[X^n/S_n]$, and let
	\[
	G_n = \langle q^n\rangle\shLog\left(1 + \sum_{i=1}^\infty F_iq^i\right).
	\]
	We then have an isomorphism of sheaves on $[X^n/S_n]$
	\[
	G_n \cong \bigoplus_{k=0}^{n-1} G_{n,k}[k]
	\]
\end{nlemma}
\begin{proof}
	By Definition \eqref{eqn:sheafLogarithm}, we have
	\[
	G_n = F_n \oplus \left( \bigoplus_{\mathfrak A} G_{\mathfrak A} \right)[1],
	\]
	where the sum runs over all set partitions $\mathfrak A$ of $[n]$ of length at least 2.
	First note that $G_{n,0} = F_n$, matching the first term.
	
	By induction, we know that for $m < n$, we have
	\begin{equation}
		\label{eqn:inductionInLogFormula}
	G_m \cong \bigoplus_{k= 0}^{m-1} G_{m,k}[k].
	\end{equation}
	Given a set partition $\mathfrak A = \{A_1,\ldots, A_l\}$ of $[n]$, with the $A_i$ forming an increasing sequence in the binary ordering, we have an isomorphism
	\begin{equation}
		\label{eqn:formulaForGMathfrakA}
	G_{\mathfrak A} = G_{A_1} \otimes \cdots \otimes G_{A_l}.
	\end{equation}
	Here each factor $G_{A_i}$ is a pullback of $G_{|A_i|}$, and so replacing $G_{|A_i|}$ by the right hand side of \eqref{eqn:inductionInLogFormula}, and further replacing each factor $G_{|A_i|,k}$ by the right hand side of \eqref{eqn:definitionOfGnk}, we express $G_{A_i}$ as a direct sum over $\mathcal T_{|A_i|}$, where each summand is a pullback of $F_T$ for $T \in \mathcal T_{|A_i|}$.
	Inserting in \eqref{eqn:formulaForGMathfrakA} we get an expression of $G_{\mathfrak A}$ as a direct sum indexed over $\mathcal T_{|A_1|} \times \cdots \times \mathcal T_{|A_l|}$.
	
	Explicitly, for an $l$-tuple $(T_i) \in \prod_{i=1}^l \mathcal T_{|A_i|}$, the corresponding summand of $G_{\mathfrak A}$ is given by
	\[
	G_{\mathfrak A, T_1,\ldots, T_l} = (F_{T_1})_{A_1}[k_1] \otimes \cdots \otimes (F_{T_l})_{A_l}[k_l],
	\]
	where $k_i$ denotes the number of non-leaf nodes of $T_i$.
	We now define a tree $T \in \mathcal T_n$ by ``gluing in the trees'' $T_i$ under the root node.
	More precisely, define $T$ such that the root node has $l$ children $v_1,\ldots, v_l$, with $L_{v_i} = A_i$ for all $i$, and such that the tree of descendants of node $v_i$ is equal to the index tree $T_i$ under the order-preserving bijection $A_i \leftrightarrow [|A_i|]$.
	Note that $T$ then has $k = 1 + \sum_{i=1}^lk_i$ non-leaf nodes.
	There is a natural isomorphism
	\begin{align*}
\gamma \colon	G_{\mathfrak A, T_1,\ldots, T_l}[1] &= ((F_{T_1})_{A_1}[k_1] \otimes \cdots \otimes (F_{T_l})_{A_l}[k_l])[1] = \left((F_{T_1})_{A_1} \otimes \cdots \otimes (F_{T_l})_{A_l}\right)[k] \\ 
=& (F_{\mathfrak A(T_1)})_{A_1} \otimes \cdots \otimes (F_{\mathfrak A(T_l)})_{A_l}[k] \cong F_{\cup \mathfrak A(T_i)}[k] = F_{\mathfrak A(T)}[k] = F_T[k],
	\end{align*}
	where the $\cong$ isomorphism uses the commutativity of the tensor product.
	
	It is easy to see that the data $(\mathfrak A, T_1,\ldots, T_l)$ and $T \in \mathcal T_n$ are in bijection by the above procedure, and so summing the isomorphisms $\gamma$ gives an isomorphism of sheaves $G_n \to \bigoplus_{k=0}^{n-1} G_{n,k}[k]$ on $X^n$.
	
	We now check that this isomorphism is $S_n$-invariant.
	For this we choose $\mathfrak A, T_1,\ldots, T_l$ and let $T$ be determined by this as above.
	Let $\sigma \in S_n$.
	For $i = 1, \ldots, l$, there is a unique $\tau_i \in S_{|A_i|}$ such that
	\[
	\sigma^{-1} \circ i_{A_i} = i_{\sigma^{-1}(A_i)} \circ \tau_i^{-1},
	\]
	where $i_{A_i} \colon [|A_i|] \to [n]$ is the order-preserving inclusion with image $A_i$.
	The $S_n$-equivariant structure on $G_n$ is given by the sum over all tuples $(\mathfrak A, T_1,\dots, T_l)$ of the isomorphism
	\begin{align*}
	\rho_\sigma \colon (\sigma_X^{-1})^*G_{\mathfrak A, T_1,\ldots, T_l} &= (\sigma_X^{-1})^*F_{T_1,A_1}[k_1] \otimes \cdots \otimes (\sigma_X^{-1})^*F_{T_l,A_l}[k_l] \\
	 &= ((\tau_1^{-1})^*F_{T_1})_{\sigma^{-1}A_1}[k_1] \otimes \cdots \otimes ((\tau_l^{-1})^*F_{T_l})_{\sigma^{-1}A_l}[k_l] \\
	 &\cong (F_{{\tau_1}^{-1}T_1})_{\sigma^{-1}A_1}[k_1] \otimes \cdots \otimes  (F_{{\tau_l}^{-1}T_l})_{\sigma^{-1}A_l}[k_l]\\
	  &= (F_{{\tau_{i_1}}^{-1}T_{i_1}})_{\sigma^{-1}A_{i_1}}[k_{i_1}] \otimes \cdots \otimes (F_{{\tau_{i_l}}^{-1}T_{i_l}})_{\sigma^{-1}A_{i_l}}[k_{i_l}] \\
	 & = G_{\sigma^{-1}(\mathfrak A), \tau_{i_1}^{-1}T_{i_1},\ldots, \tau^{-1}_{i_l}T_{i_l}}
	\end{align*}
	where $\sigma^{-1}A_{i_1},\ldots, \sigma^{-1}A_{i_l}$ are increasing in the binary ordering, and where we use $\rho_{T_i,\tau_i}$ from \eqref{eqn:rhoSigma} in the third isomorphism.

	We now have a diagram of isomorphisms
	\[
	\begin{tikzcd}
		(\sigma_X^{-1})^*G_{\mathfrak A, T_1,\ldots, T_l}[1] \arrow[r, "\rho_\sigma"] \arrow[d, "\gamma"] & G_{\sigma^{-1}(\mathfrak A), \tau_{i_1}^{-1}(T_{i_1}),\ldots, \tau^{-1}_{i_l}(T_{i_l})}[1] \arrow[d, "\gamma"] \\
		(\sigma_X^{-1})^*F_T[k] \arrow[r, "\rho_{T,\sigma}"] & F_{\sigma^{-1}(T)}[k],
	\end{tikzcd}
	\]
	which we claim is commutative.

	Let $\mathfrak A(T) = \{B_1,\ldots, B_j\}$.
	Up to a sign, both compositions in the diagram equal the isomorphism obtained by
	\begin{align*}
	(\sigma_X^{-1})^*G_{\mathfrak A, T_1,\ldots, T_l} &= ((\sigma_X^{-1})^*F_{T_1,A_1}[k_1] \otimes \ldots \otimes (\sigma_X^{-1})^*F_{T_l,A_l}[k_l])[1] \\
	&= ((\sigma_X^{-1})^*F_{T_1,A_1} \otimes \ldots \otimes (\sigma_X^{-1})^*F_{T_l,A_l})[k] \\
	&= (\sigma_X^{-1})^*(F_{B_1} \otimes \cdots \otimes F_{B_j})[k] \\
	&\cong (F_{\sigma^{-1} B_1} \otimes \cdots \otimes F_{\sigma^{-1}B_j}) \\
	&= (F_{\sigma^{-1}B_{l_1}} \otimes \cdots \otimes F_{\sigma^{-1}B_{l_j}})[k] = F_{\sigma^{-1}(T)}[k].
	\end{align*}
	Here the $l_i$ are such that $\sigma^{-1}B_{l_i}$ are increasing in the binary ordering.
	The fourth isomorphism is the tensor product of the $\rho_{B_i,\sigma}$, while the other ones are the ``standard'' isomorphisms, including the commutation isomorphism of the tensor product.

	Going around the lower left corner of the diagram, the definition \eqref{eqn:rhoSigma} of $\rho_{T,\sigma}$ means that the above isomorphism is modified by a sign $(-1)^{l(T,\sigma)}$ where $l(T,\sigma)$ is the number of pairs $A < B \in P(T)$ such that $\sigma(A) > \sigma(B)$.
	Going around the upper right corner, the above isomorphism must be modified in two ways, first a sign $(-1)^{\sum l(T_i,\tau_i)}$, secondly, because of Remark \ref{remark:noncommutativity}, a sign $(-1)^s$, where 
	\[
	s = \sum k_jk_{j'},
	\]
	the sum running over all $1 \le j < j' \le l$ such that $i_{j'} < i_j$.
	Now since $k_j = |P(T_j)|$, one can check that 
	\[
	s + \sum l(T_i,\tau_i) = l(T,\sigma).
	\]
Thus the sign modifications around both corners of the diagram agree, and the diagram commutes.
\end{proof}	

\section{Systems of sheaves}
\label{sec:systems}
If $(F_n)_{n=1}^\infty$ is a factorisable sequence of sheaves on $[X^n/S_n]$, then for any given $n$ and set partitions $\mathfrak A, \mathfrak B$ of $[n]$, there is an isomorphism $F_{\mathfrak A} \to F_{\mathfrak B}$ on some open $U_{\mathfrak A, \mathfrak B} \subseteq X^n$ (see Lemma \ref{thm:factorisableGivesSystem}).
The goal of this section is to prove Proposition \ref{thm:systemToGoodSystem}, which says that after replacing the $F_{\mathfrak A}$ by K-theoretically equivalent sheaves $F'_{\mathfrak A}$, we may assume that these isomorphisms can be extended to homomorphisms defined on the whole of $X^n$ (preserving compatibilities under composition and the $S_n$-action).

Because this result is also applied in a more general setting in \cite{thimm_ktheoretic_2024}, we abstract the setup.
Let $Y$ be a quasi-projective scheme, carrying actions by an algebraic group $G$ and a finite group $S$, which we assume to commute.
In this section, we write ``sheaf'' for $G$-equivariant, cohomologically $\ZZ/2$-graded coherent sheaf.

Our main example has $Y = X^n, G = \mathbb{T}$ and $S = S_n$.
\subsection{The stratification and set of sheaves}
Let $\mathcal{I}$ be a finite set, and for each $\alpha \in \mathcal{I}$, let $Y_{\alpha} \subseteq Y$ be a $G$-invariant locally closed subset of $Y$.
We require the collection of $Y_\alpha$ to form a stratification of $Y$, meaning they satisfy the following conditions:
\begin{itemize}
	\item Every point of $Y$ lies in a unique $Y_\alpha$.
	\item The closure of every $Y_\alpha$ is a union of $Y_\beta$.
\end{itemize}
We define a partial ordering on $\mathcal{I}$ by letting $\alpha \ge \beta$ if $Y_\alpha \subseteq \overline{Y}_\beta$.

We further assume that the $S$-action preserves the collection of $Y_\alpha$, and so we can define an action of $S$ on $\mathcal{I}$ by
\[
\sigma(Y_\alpha) = Y_{\sigma\alpha}.
\]
Define the open subset $U_\alpha \subseteq Y$ by
\[
U_\alpha = \bigcup_{\beta \le \alpha} Y_\beta.
\]

Let now $\mathcal{J}$ be a finite partially ordered set, equipped with the following data.
\begin{itemize}
	\item An action of $S$ on $\mathcal{J}$.
	\item For every $\alpha \in \mathcal{I}$, an equivalence relation $i \sim_\alpha j$ on $\mathcal{J}$.
\end{itemize}
These satisfy the following axioms:
\begin{itemize}
	\item[(A1)] The action of $S$ on $\mathcal{J}$ preserves the partial ordering.
	\item[(A2)] For every $\sigma \in S$, $i, j \in \mathcal{J}$ and $\alpha \in \mathcal{I}$, we have
	\[
	i \sim_{\alpha} j \Leftrightarrow \sigma i \sim_{\sigma\alpha} \sigma j,
	\]
	\item[(A3)] If $\alpha \ge \beta$, then $i \sim_\alpha j \Rightarrow i \sim_{\beta} j$.
	\item[(A4)] If $i \le j \le k$, then $i \sim_\alpha k$ is equivalent to $i \sim_\alpha j$ and $j \sim_\alpha k$.
	\item[(A5)] Let $\alpha, \beta \in \mathcal{I}$.
	Given elements $i, j, k \in \mathcal{J}$ with the relations as in the two leftmost diagrams, there exists a $k' \in \mathcal{J}$ with relations as in the two rightmost diagrams.
	\[
	\begin{tikzcd}
 & i \arrow[d, phantom, "\sim_\beta"{rotate = -90}] &		& i \arrow[d, phantom, "\le"{rotate = -90}] & k' \arrow[d, phantom, "\sim_\beta"{rotate = -90}] \arrow[r, phantom, "\sim_\alpha"] & i \arrow[d, phantom, "\sim_\beta"{rotate = -90}] &	k'\arrow[r, phantom, "\le"] \arrow[d, phantom, "\le"{rotate = -90}]	& i \arrow[d, phantom, "\le"{rotate = -90}]  \\
k \arrow[r, phantom, "\sim_\alpha"]& j &	k \arrow[r, phantom, "\leq"]	& j & k \arrow[r, phantom, "\sim_\alpha"]& j &	k \arrow[r, phantom, "\leq"]	& j
	\end{tikzcd}
	\]
\end{itemize}

\begin{nex}
	\label{ex:mainExampleSystemContext}
	In our main example, we have $Y = X^n$, $S = S_n$ and $G = \mathbb{T}$, with the $S_n$- and $\mathbb{T}$-actions on $X^n$ the usual ones.
	We let $\mathcal{I}$ and $\mathcal{J}$ both be the set of all set partitions of $[n]$, equipped with the standard action of $S_n$, and with the partial ordering for which $\mathfrak A \le \mathfrak B$ if $\mathfrak A$ is finer than $\mathfrak B$, that is if every element of $\mathfrak A$ is contained in some element of $\mathfrak B$.
	
	For $\mathfrak A \in \mathcal{I}$, we let $(X^n)_{\mathfrak A}$ be the set of $(x_1,\ldots, x_n)$ such that $x_i = x_j$ if and only if $i$ and $j$ are contained in the same element of $\mathfrak A$.

	For two set partitions $\mathfrak A, \mathfrak B$ of $[n]$, we write $\mathfrak A \cap \mathfrak B$ for the set partition whose elements are all the non-empty intersections of elements in $\mathfrak A$ and $\mathfrak B$.
	Equivalently, $\mathfrak A \cap \mathfrak B$ is the coarsest set partition finer than both $\mathfrak A$ and $\mathfrak B$.

	Given $\mathfrak A \in \mathcal{I}$, we define the equivalence relation $\sim_{\mathfrak A}$ on $\mathcal{J}$ by, for $\mathfrak B, \mathfrak C \in \mathcal{J}$, saying
	\[
	\mathfrak B \sim_{\mathfrak A} \mathfrak C \overset{\text{def}}{\Leftrightarrow} \mathfrak B \cap \mathfrak A = \mathfrak C \cap \mathfrak A.
	\]
	
	Axioms (A1), (A2), (A3), (A4) are now all easy to verify.
	For axiom (A5), letting $\alpha, \beta \in \mathcal{I}$ and $i,j,k \in \mathcal{J}$ be set partitions of $[n]$ satisfying the requirements, we can take $k' = i \cap k$, and then we get
	\[
	k' \cap \beta = i \cap k \cap \beta = j \cap k \cap \beta = k \cap \beta,
	\]
	using $i \sim_\beta j$ and $k \le j$.
	We also get 
	\[
	k' \cap \alpha = i \cap k \cap \alpha = i \cap j \cap \alpha = i \cap \alpha,
	\]
	using $k \sim_\alpha j$ and $i \le j$.
	This shows that $k' \sim_\beta k$ and $k' \sim_\alpha i$, so (A5) is satisfied.
\end{nex}

\subsection{Systems of sheaves}
For $i, j \in \mathcal{J}$, we let 
\[
U_{i,j} = \bigcup_{\substack{\alpha \in \mathcal{I} \\ i \sim_\alpha j}} U_{\alpha} = \bigcup_{\substack{\alpha \in \mathcal{I} \\ i \sim_\alpha j}} Y_{\alpha},
\]
the second equality by (A3).

\begin{ndefn}
	A \textbf{system} (of sheaves) with respect to the above is defined by giving the following data:
	\begin{itemize}
		\item For each $i \in \mathcal{J}$, a sheaf $F_i$ on $Y$.
		\item For each $i \le j \in \mathcal{J}$ an isomorphism of sheaves
		\[
		\phi_{i,j} \colon F_i|_{U_{i,j}} \to F_j|_{U_{i,j}}.
		\]
	\end{itemize}
	We denote a system by $(F,\phi)$.
	We require the following two conditions
	\begin{itemize}
		\item[(C1)] For any $i \in \mathcal{J}$, we have $\phi_{i,i} = \id_{F_i}$. 
		\item[(C2)] For any $i \le j \le k \in \mathcal{J}$, over $U_{i,j} \cap U_{j,k}$,\footnote{Note that by assumption (A4), we have $U_{i,k} = U_{i,j} \cap U_{j,k}$.} we have
		\[
		\phi_{j,k} \circ \phi_{i,j} = \phi_{i,k}.
		\]
	\end{itemize}
\end{ndefn}

A homomorphism between two systems of sheaves $(F_i, \phi_{i,j})$ and $(G_i,\psi_{i,j})$ is given by a homomorphism from each $F_i$ to each $G_i$ commuting with the isomorphisms $\phi$ and $\psi$.
This gives an abelian category $\textbf{Syst}$\footnote{The category $\Syst$ depends on the $G$-equivariant scheme $Y$, the posets $\mathcal{I}, \mathcal{J}$, the strata $\{Y_\alpha\}_{\alpha \in \mathcal{I}}$ and the equivalence relations $\{\sim_\alpha\}_{\alpha \in \mathcal{I}}$ on $\mathcal{J}$}, where a sequence
\[
(F',\phi') \to (F, \phi) \to (F'',\phi'')
\]
is exact if and only if the underlying sequence of sheaves
\[
F'_i \to F_i \to F''_i
\]
is exact for all $i \in \mathcal{J}$.

The group $S$ acts on $\mathbf{Syst}$ by, for $\sigma \in S$, letting
\[
\sigma(F, \phi) = (\sigma(F), \sigma(\phi))
\]
where
\[
\sigma(F)_i = (\sigma^{-1})^*F_{\sigma^{-1}(i)},
\]
and
\[
\sigma(\phi)_{i,j} = (\sigma^{-1})^*(\phi_{\sigma^{-1}(i),\sigma^{-1}(j)}).
\]
\begin{ndefn}
	An $S$-\textbf{equivariant system} of sheaves is an object in the equivariant category $\mathbf{Syst}^S$ (see e.g.~ \cite[Def.~3.1.]{beckmann_equivariant_2023}).
	Explicitly, this is given by an $(F,\phi) \in \Syst$, and for every $\sigma \in S, i \in \mathcal{J}$, an isomorphism
\[
\rho_{\sigma, i} \colon (\sigma^{-1})^*F_{\sigma^{-1}(i)} \to F_i
\]
such that
\begin{itemize}
	\item for all $i \le j \in \mathcal{J}$ and $\sigma \in S$, we have $\rho_{\sigma, j} \circ \phi_{\sigma^{-1}(i), \sigma^{-1}(j)} = \phi_{i,j} \circ \rho_{\sigma,i}$.
	\item for all $\sigma, \psi \in S$ and $i \in \mathcal{J}$, the isomorphism
\[
((\sigma\psi)^{-1})^*(F_{(\sigma\psi)^{-1}(i)}) = (\psi^{-1})^*(\sigma^{-1})^*F_{\psi^{-1}(\sigma^{-1}(i))} \to (\sigma^{-1})^*F_{\sigma^{-1}(i)} \to F_i
\]
equals $\rho_{\sigma\psi, i}$.
\end{itemize}
\end{ndefn}

\begin{nex}
	As we will show in Lemma \ref{thm:factorisableGivesSystem}, if $(F_n)_{n=1}^\infty$ is a factorisable sequence of sheaves on $[X^n/S_n]$, then the collection of sheaves $F_{\mathfrak A}$ for $\mathfrak A$ a set partition of $[n]$ has the structure of an $S_n$-equivariant system on $X^n$.
\end{nex}

\subsection{Strict systems of sheaves}
\begin{ndefn}
	A \textbf{strict system} (of sheaves) is defined by the following data:
	\begin{itemize}
		\item For each $i \in \mathcal{J}$, a sheaf $F_i$ on $Y$.
		\item For each $i \le j \in \mathcal{J}$ a homomorphism $\phi_{i,j} \colon F_i \to F_j$.
	\end{itemize}
	This data must satisfy conditions:
	\begin{itemize}
		\item[(D1)] For any $i \le j \le k \in \mathcal{J}$, we have $\phi_{j,k} \circ \phi_{i,j} = \phi_{i,k}$.
		\item[(D2)] $\phi_{i,i}$ is the identity for all $i \in \mathcal{J}$.
		\item[(D3)] If $i \le j$, then $\phi_{i,j}|_{U_{i,j}}$ is an isomorphism.
	\end{itemize}
\end{ndefn}
More concisely, thinking of the poset $\mathcal{J}$ as a category, a strict system is a functor $\mathcal{J} \to \Coh^G_{\ZZ/2}(Y)$ satisfying the extra condition given by (D3).

In the same way as for systems, there is a natural abelian category $\SSyst$ whose objects are strict systems, and there is an exact functor $\Psi \colon \SSyst \to \Syst$ given by
\[
\Psi(F_i, \phi_{i,j}) = (F_i, \phi_{i,j}|_{U_{i,j}}).
\]

\subsection{The functor $D_\alpha$}
We would like to show that in $K$-theory, the functor $\Psi$ (and its $S$-equivariant version) induces a surjection.
The idea is to produce, for any system $(F,\phi)$, a strict system $(F',\phi')$ and a morphism $(F,\phi) \to \Psi(F',\phi')$ whose kernel and cokernel are ``smaller'' than $(F,\phi)$; inducting on the ``size'' of $(F,\phi)$ then shows that the K-theory class of $(F,\phi)$ is in the image of $\Psi$.

A first approximation to the construction of $(F',\phi')$ goes as follows.  Let $\alpha \in \mathcal{I}$, and let $j_\alpha \colon Y_{\alpha} \to Y$ be the inclusion.
Given a system of sheaves $(F, \phi)$, we could try to take $F'_i = (j_{\alpha})_*j_\alpha^*F_i$ and
\[
\phi'_{i,j} = \begin{cases} (j_\alpha)_*j_\alpha^*\phi_{i,j} \text{ if }i \sim_\alpha j \\
	0 \text{ if } i \not\sim_\alpha j.
	\end{cases}
\]
This almost defines a strict system, the only problem being that $(j_\alpha)_*j_\alpha^*(F_i)$ may fail to be a coherent sheaf, since $Y_\alpha$ is not necessarily closed in $Y$.

We will modify the above construction to give a valid functor $D_\alpha \colon \Syst \to \SSyst$.
Roughly speaking, $D_\alpha(F,\phi)$ is the smallest strict subobject of the $(F',\phi')$ defined above for which the morphism $(F,\phi) \to \Psi(F',\phi')$ factors through $\Psi(D_\alpha(F,\phi))$.
\begin{ndefn}
	Let $\alpha \in \mathcal{I}$. 
	Define a functor $D_\alpha \colon \Syst \to \SSyst$ by, for $i \in \mathcal{J}$, setting
	\[
	D_{\alpha}(F)_i = \image\left(\bigoplus_{\substack{k \sim_\alpha i \\ k \le i}} F_k \to \bigoplus_{\substack{k \sim_\alpha i\\ k \le i}} (j_\alpha)_*j_\alpha^*F_k  \to (j_\alpha)_*j_\alpha^* F_i\right),
	\]
	where the final arrow is the sum of all the homomorphisms $(j_{\alpha})_*j^*_\alpha(\phi_{k,i})$.

For any $i \le j \in \mathcal{J}$, let 
\[
D_{\alpha}\phi_{i,j} \colon D_{\alpha}(F)_i \to D_{\alpha}(F)_j
\]
be defined by
\begin{itemize}
	\item If $i \not\sim_\alpha j$, let $D_{\alpha}\phi_{i,j} = 0$.
	\item If $i \sim_\alpha j$, then the commutative diagram
	\[
	\begin{tikzcd}
		\bigoplus\limits_{\substack{k \sim_\alpha i \\ k \le i}} F_k \arrow[r] \arrow[d] & \bigoplus\limits_{\substack{k \sim_\alpha i\\k \le i}}(j_\alpha)_*j^*_\alpha F_k \arrow[r] \arrow[d] & (j_\alpha)_*j^*_\alpha F_i \arrow[d, "(j_\alpha)_*j_\alpha^*\phi_{i,j}", "\cong"'] \\
		\bigoplus\limits_{\substack{k \sim_\alpha j \\k \le j}} F_k \arrow[r] & \bigoplus\limits_{\substack{k \sim_\alpha j \\ k \le j}} (j_\alpha)_*j^*_\alpha F_k \arrow[r] & (j_\alpha)_*j^{*}_\alpha F_j
	\end{tikzcd}
	\]
	induces a homomorphism between the images of the composition in the top and bottom row, which defines $D_{\alpha}\phi_{i,j}$.
\end{itemize}
\end{ndefn}

To verify that the above definition actually defines an object of $\SSyst$, note that (A4) easily implies that (D1) holds, and (D2) trivially holds.
The fact that (D3) holds is the following lemma.
\begin{nlemma}
	Let $(F,\phi) \in \Syst$, let $\alpha \in \mathcal{I}$ and $i \le j \in \mathcal{J}$ be such that $i \sim_\alpha j$.
	Then $D_\alpha\phi_{i,j}|_{U_{i,j}}$ is an isomorphism.
\end{nlemma}
\begin{proof}
It is enough to show that $D_\alpha\phi_{i,j}|_{U_{\beta}}$ is an isomorphism for each $\beta \in \mathcal{I}$ with $i \sim_\beta j$, since $U_{i,j}$ is by definition the union of these $U_\beta$.

It is clear that $D_{\alpha} F_i$ and $D_{\alpha} F_j$ are both supported on $\overline Y_\alpha$.
If $\beta \not \ge \alpha$, then $U_\beta \cap \overline Y_\alpha = \varnothing$, and so
\[
D_{\alpha}F_i|_{U_\beta} = D_{\alpha}F_j|_{U_\beta} = 0,
\]
and so $D_{\alpha}\phi_{i,j}|_{U_\beta}$ is an isomorphism.
	
We now assume that $\beta \ge \alpha$.
	The homomorphism
	\[
	\bigoplus_{k \sim_\alpha i, k \le i} F_k \to \bigoplus_{k \sim_\alpha j, k \le j}F_k
	\]
	is an inclusion of summands, hence injective, and it follows that $D_{\alpha}\phi_{i,j}$ is injective as well.
	To show surjectivity over $U_\beta$, it is enough to show that for every $F_k$ in the second direct sum there exists an $F_{k'}$ in the first direct sum whose image in $(j_\alpha)_*j_\alpha^*F_j$ agrees with that of $F_k$ over $U_\beta$.
	By axiom (A5), there exists a $k'$ such that $k' \le k$, $k' \sim_\beta k$, $k' \le i$ and $k' \sim_\alpha i$.
We then get a commutative diagram
	\[
	\begin{tikzcd}
		F_{k'}|_{U_\beta} \arrow[d] \arrow[r] & ((j_\alpha)_*j_\alpha^*F_{k'})|_{U_\beta} \arrow[r] \arrow[d] & ((j_\alpha)_*j_\alpha^*F_i)|_{U_\beta} \arrow[d] \\
		F_k|_{U_\beta} \arrow[r] & ((j_\alpha)_*j_\alpha^*F_k)|_{U_\beta} \arrow[r] & ((j_\alpha)_*j_\alpha^* F_j)|_{U_\beta},
	\end{tikzcd}
	\]
	It follows that over $U_\beta$, the image of $F_k$ in $(j_\alpha)_*j_\alpha^*F_j$ is the same as that of $F_{k'}$, so $\phi_{i,j}|_{U_\beta}$ is surjective.
\end{proof}

For each $i \in \mathcal{J}$, there is a homomorphism $I_{\alpha,i} \colon F_i \to D_{\alpha}(F)_i$ obtained from the inclusion of $F_i$ in $\oplus_{k \le i, k \sim_\alpha i}F_k$.
One checks that this gives a homomorphism of systems
\[
I_\alpha \colon (F,\phi) \to \Psi D_\alpha(F,\phi).
\]

\subsection{The functor $\widetilde{D}_\alpha$}
The group $S$ acts on the category $\SSyst$ in the same way as on $\Syst$.
An $S$-equivariant strict system is an $S$-equivariant object in $\SSyst$, that is an object of $\SSyst^S$.

For any $(F,\phi) \in \SSyst$, $\beta \in \mathcal{I}$ and $\sigma \in S$, there is a natural isomorphism
\[
\sigma(D_\beta(F,\phi)) \cong D_{\sigma \beta}(\sigma(F,\phi)).
\]
If $(F,\phi) \in \SSyst^S$, the structure isomorphism $\sigma(F,\phi) \cong (F,\phi)$ composes with the above to define an isomorphism
\begin{equation}
	\label{eqn:DBetaEquivariance}
\sigma(D_\beta(F,\phi)) \cong D_{\sigma\beta}(F,\phi).
\end{equation}

We now modify the functor $D_\alpha$ to a functor
\[
\widetilde D_\alpha \colon \Syst^S \to \SSyst^S,
\]
defined as follows.
Let $(F,\phi) \in \Syst^S$, and let
\[
\widetilde{D}_\alpha(F,\phi) = \bigoplus_{\beta \in S\alpha} D_\beta(F,\phi),
\]
equipped with an $S$-equivariant structure by taking the sum of the isomorphisms \eqref{eqn:DBetaEquivariance}.
We also get a homomorphism in $\Syst^S$
\[
\widetilde{I}_\alpha = \sum_{\beta \in S\alpha} I_\beta \colon (F,\phi) \to \bigoplus_{\beta \in S\alpha} D_\beta(F,\phi) = \Psi\widetilde{D}_\alpha(F,\phi).
\]
We thus get the following exact sequence in $\mathbf{Syst}^S$
\[
0 \to \ker \widetilde{I}_\alpha \to (F,\phi) \overset{\widetilde{I}_\alpha}{\to} \Psi\widetilde{D}_\alpha(F,\phi) \to \cok \widetilde{I}_\alpha \to 0.
\]

\subsection{Noetherian induction}
We now show that $\ker \widetilde{I}_\alpha$ and $\cok \widetilde I_\alpha$ have smaller support than $F$, in a precise sense.
For $l \ge 1$, let the $l$-fold thickening of $\overline{Y}_\alpha$ be the scheme defined by the ideal sheaf $I_{\overline{Y}_\alpha}^l$.

For any $F \in \Syst$,\footnote{Since the isomorphisms $\phi_{i,j}$ play no role from this point on, we omit them from the notation.} let $l(\alpha, F) \in \mathbb N \cup \{\infty\}$ be the minimal number $l$ such that, for every $i \in \mathcal{J}$, the sheaf $F_i|_{U_\alpha}$ is scheme-theoretically supported on the $l$-fold thickening of $\overline{Y}_\alpha$.
We let $\Supp(F)$ be the union of the scheme-theoretic supports of all $F_i$.
The function $l(\alpha, F)$ has the following easily verified properties:
\begin{itemize}
	\item $l(\alpha,F) = 0$ if and only if $F|_{U_\alpha} = 0$.
	\item $l(\alpha,F) = 0$ for all $\alpha \in \mathcal{I}$ if and only if $F = 0$.
	\item If $\alpha$ is minimal in the set
	\[
	\{\alpha \in \mathcal{I} \mid l(\alpha,F) \not= 0\},
	\]
	then $0 < l(\alpha,F) < \infty$.
\end{itemize}

\begin{nlemma}
	\label{thm:inequalitiesForL}
	Let $\alpha \in \mathcal{I}$ and $F \in \Syst^S$.
	\begin{enumerate}
		\item If $0 < l(\alpha, F) < \infty$, then $l(\alpha, \ker \widetilde I_\alpha) < l(\alpha, F)$.
		\item $l(\alpha, \cok \widetilde I_\alpha) = 0$.
	\end{enumerate}
\end{nlemma}
\begin{proof}
(1): Restricting to $U_\alpha$, the homomorphism $I_\alpha \colon F \to D_{\alpha}(F)$ is for each $i \in \mathcal{J}$ naturally identified with the surjection
	\[
	F_i \to F_i|_{Y_{\alpha}},
	\]
	so that in particular we have $(\ker I_\alpha)_i = \mathcal{I}_{Y_\alpha}F_i$.
	Since $\ker \widetilde I_\alpha \subseteq \ker I_\alpha$, we then find that
	\[
	l(\alpha, \ker \widetilde{I}_\alpha) \le l(\alpha, \ker I_\alpha) = l(\alpha, \mathcal I_{Y_\alpha}F) = l(\alpha, F) - 1 < l(\alpha, F).
	\]
	
	(2): For any $\beta \in \mathcal{I}$, we have $\Supp(D_\beta(F)) \subseteq \overline Y_\beta$.
	For any $\beta \not= \alpha$ in the same $S$-orbit as $\alpha$, we have that $\overline{Y}_\beta \cap U_\alpha = \varnothing$, so that $D_\beta(F)|_{U_\alpha} = 0$.
	On $U_\alpha$, then, we may naturally identify $I_\alpha$ and $\widetilde I_\alpha$, so that 
	\[
	\cok \widetilde{I}_\alpha|_{U_\alpha} = \cok I_{\alpha}|_{U_\alpha} = 0.
	\]
\end{proof}

\begin{ncor}
	\label{thm:supportShrinks}
	Let $0 \not= F \in \Syst^S$, and let $\alpha \in \mathcal{I}$ be such that $0 < l(\alpha, F) < \infty$.
	Then 
	\[
	\Supp(\ker \widetilde I_\alpha), \Supp(\cok \widetilde I_\alpha) \subsetneq \Supp(F).
	\]
\end{ncor}
\begin{proof}
	Since $\ker \widetilde I_\alpha$ is a subobject of $F$, it is clear that $\Supp(\ker \widetilde I_\alpha) \subseteq \Supp(F)$.
	Part (1) of Lemma \ref{thm:inequalitiesForL} shows the inclusion must be strict.
	
	For $\cok \widetilde I_\alpha$, note that in general we have $\Supp(\widetilde D_\alpha(F)) \subseteq \Supp(F)$, and $\cok \widetilde I_\alpha$ being a quotient of $\widetilde D_\alpha(F)$ then gives $\Supp(\cok \widetilde I_\alpha) \subseteq \Supp(F)$.
	Part (2) of Lemma \ref{thm:inequalitiesForL} gives $0 = l(\alpha, \cok \widetilde I_\alpha) < l(\alpha, F)$, so the inclusion must be strict.
\end{proof}

\begin{ncor}
	\label{thm:KTheoryMapSurjective}
	The homomorphism 
	\[
	K(\Psi) \colon K(\mathbf{SSyst}^S) \to K(\Syst^S)
	\]
	is surjective.
\end{ncor}
\begin{proof}
	Suppose not, then we can find an $F \in \mathbf{Syst}^S$ such that $[F]$ is not in the image of $K(\Psi)$.
	Since $Y$ is Noetherian, we can find an $F$ such that $\Supp(F)$ is minimal with respect to this property.
	
	Taking $\alpha \in \mathcal{I}$ to be minimal with $l(\alpha, F) \not= 0$, we have $0 < l(\alpha,F) < \infty$.
	We can write
	\[
	[F] = [\ker \widetilde I_{\alpha}] + [\Psi\widetilde{D}_\alpha(F)] - [\cok \widetilde I_\alpha].
	\]
	The minimality property of $F$ and Corollary \ref{thm:supportShrinks} together imply that all the terms on the right hand side are in the image of $K(\Psi)$, which gives a contradiction.
\end{proof}

\section{From factorisable sequences to systems of sheaves}
\begin{ndefn}
	\label{def:factorisableSequence}
	Let $(\cF_n)_{n=1}^{\infty}$ be a sequence where each $\cF_n$ is a ($\mathbb T$-equivariant, $\ZZ/2$-graded) coherent sheaf on $\Sym^n(X)$. Then \emph{factorisation data} for $(\cF_n)_{n=1}^{\infty}$ consists of a collection of $T$-equivariant isomorphisms
	$$
	\phi_{m,n} \colon (\cF_m \boxtimes \cF_n)|_{U_{m,n}} \stackrel{\cong}{\to} \cF_{m+n}|_{U_{m,n}},
	$$
	for all $m,n \geq 1$, satisfying the following properties:
	\begin{enumerate}
		\item \textbf{(Associativity)} For any $m,n,p > 0$, we have $$\phi_{m,n+p} \circ (\id_{\cF_m} \boxtimes \phi_{n,p}) = \phi_{m+n,p} \circ (\phi_{m,n} \boxtimes \id_{\cF_p})$$ as isomorphisms on the set 
		\[
		U_{m,n,p} \subseteq \Sym^m X \times \Sym^n X \times \Sym^p X
		\]
		where the three 0-cycles have pairwise disjoint support.
		\item \textbf{(Commutativity)} Let $m,n >0$ and denote by $$\tau \colon \Sym^m(X) \times \Sym^n(X) \to \Sym^n(X) \times \Sym^m(X)$$ the permutation map. Then the following diagram commutes
		\[
		\begin{tikzcd}
			(\cF_m \boxtimes \cF_n) |_{U_{m,n}} \arrow[r,"\phi_{m,n}"] \arrow[d, "\cong", "v"'] & \cF_{m+n}|_{U_{m,n}} \arrow["\cong", "w"', d] \\ 
			\tau^*((\cF_n \boxtimes \cF_m) |_{U_{n,m}}) \arrow["\tau^*\phi_{n,m}",r] & \tau^*(\cF_{n+m}|_{U_{n,m}})
		\end{tikzcd}
		\]
		where the vertical arrows are canonical isomorphisms induced by the canonical isomorphism $\cF_m \boxtimes \cF_n \cong \tau^* (\cF_n \boxtimes \cF_m)$ and the commutative diagram
		\[
		\begin{tikzcd}
			U_{m,n} \arrow[r, "\tau"] \arrow["\sigma_{m,n}"',dr] & U_{n,m} \arrow[d, "\sigma_{n,m}"] \\
			& \Sym^{m+n}(X)
		\end{tikzcd}
		\]
	\end{enumerate}
	We say that $(\cF_n)_{n=1}^{\infty}$ is \emph{factorisable} if there exists factorisation data for $(\cF_n)_{n=1}^{\infty}$.
\end{ndefn}

\begin{ndefn}
	Let $(F_n)_{n=1}^\infty$ be a sequence where each $F_n$ is a sheaf on $[X^n/S_n]$.
	We say that this sequence is factorisable if there exists a factorisable sequence $(\cF_n)_{n=1}^\infty$ in the sense of Definition \ref{def:factorisableSequence} such that for each $n$ we have $F_n = p^{*}(\cF_n)$ with $p \colon [X^n/S_n] \to \Sym^n X$ the coarse moduli space morphism.
\end{ndefn}

In the following lemma, a system means an object of the category defined by Example \ref{ex:mainExampleSystemContext}.
\begin{nlemma}
	\label{thm:factorisableGivesSystem}
	Let $(F_n)_{n=1}^\infty$ be a factorisable sequence of sheaves on $[X^n/S_n]$.
	For every $n$, there exists an $S_n$-equivariant system on $X^n$ where the sheaves are $F_{\mathfrak A}$, and the $S_n$-structure is the one induced by $\rho_{\sigma, \mathfrak A}$.
\end{nlemma}
\begin{proof}
	For any set partitions $\mathfrak A \le \mathfrak B$, we must define isomorphisms of sheaves
	\[
	\phi_{\mathfrak A, \mathfrak B} \colon F_{\mathfrak A}|_{U_{\mathfrak A, \mathfrak B}} \to F_{\mathfrak B}|_{U_{\mathfrak A, \mathfrak B}},
	\]
	satisfying conditions (C1), (C2), as well as the $S_n$-equivariance compatibility of, for each $\sigma \in S_n$, the commutativity of the following diagram.
	\[
	\begin{tikzcd}
		(\sigma^{-1}_X)^*F_{\mathfrak A} \arrow[r, "\rho_{\sigma,\mathfrak A}"] \arrow[d, "(\sigma_X^{-1})^*\phi_{\mathfrak A, \mathfrak B}"] & F_{\sigma^{-1}\mathfrak A} \arrow[d, "\phi_{\sigma^{-1}\mathfrak A, \sigma^{-1}\mathfrak B}"] \\
		(\sigma_X^{-1})^*F_{\mathfrak B} \arrow[r, "\rho_{\sigma,\mathfrak B}"] & F_{\sigma^{-1}\mathfrak B}
	\end{tikzcd}
	\]
	
	Since $(F_n)_{n=1}^\infty$ is a factorisable sequence, there exist sheaves $\cF_n$ on $\Sym^n X$ and partial isomorphisms $\phi_{m,n}$ as in Def.~\ref{def:factorisableSequence} such that $F_n = p^*\cF_n$ for $p \colon [X^n/S_n] \to \Sym^n X$.
	For any two disjoint subsets ${A_1}, {A_2} \subseteq [n]$, we have the commutative diagram
	\[
	\begin{tikzcd}
		X^n \arrow[r] & X^{A_1} \times X^{A_2} \arrow[r, "\cong"] \arrow[d] & X^{{A_1} \cup {A_2}} \arrow[d] \\
		& \Sym^{|{A_1}|}X \times \Sym^{|{A_2}|} X \arrow[r] & \Sym^{|{A_1}| + |{A_2}|}X.
	\end{tikzcd}
	\]
	Pulling back the isomorphism $\phi_{|{A_1}|, |{A_2}|}$ to $X^n$ gives an isomorphism
	\[
	\phi_{\{{A_1},{A_2}\}} \colon F_{A_1} \otimes F_{A_2} \to F_{{A_1} \cup {A_2}},
	\]
	defined on 
	\[
	U_{A_1,A_2} := \{(x_i)_{i=1}^n \mid j \in A_1, k \in A_2 \Rightarrow x_j \not= x_k\} \subseteq X^n.
	\]
	Let now $A_1,\ldots, A_l \subseteq [n]$ be disjoint.
	For each $1 \le i \le l-1$, tensoring the above isomorphism with identity isomorphisms lets us define an isomorphism
	\[
	F_{A_1 \cup \cdots \cup A_{i}} \otimes F_{A_{i+1}} \otimes \cdots \otimes F_{A_l} \to F_{A_1 \cup \cdots \cup A_{i+1}} \otimes F_{A_{i+2}} \otimes \cdots \otimes F_{A_l}
	\]
	The composition of these gives an isomorphism
	\[
	F_{A_1} \otimes \cdots \otimes F_{A_l} \to F_{A_1 \cup \cdots \cup A_l}
	\]
	defined on
	\[
	U_{A_1,\ldots, A_l} := \bigcap_{i \not= j} U_{A_i,A_j}.
	\]
	Using the commutativity and associativity axioms of \cite[Def.~2.5]{kool_proof_2025}, one checks that this isomorphism commutes with the tensor commutativity isomorphism under reordering of the $A_i$.
	We write $\phi_{\{A_1,\ldots, A_l\}}$ for the above isomorphism.
	
	For any two set partitions of $[n]$,
	\[
	\{A_{i,j}\}_{i,j} = \mathfrak A \le \mathfrak B = \{B_i\}_i,
	\]
	where $B_i = \cup_j A_{i,j}$ for all $i$, taking the tensor product of the isomorphisms $\phi_{\{A_{i,j}\}_j}$ defines an isomorphism $\phi_{\mathfrak A, \mathfrak B} \colon F_{\mathfrak A} \to F_{\mathfrak B}$ over $U_{\mathfrak A, \mathfrak B}$.

	Condition (C1) is obvious, condition (C2) is easy but tedious to verify, as is the fact that the system is $S_n$-equivariant.
	We omit the details.
\end{proof}

\begin{nprop}
	\label{thm:systemToGoodSystem}
	If $(F_n)_{n=1}^\infty$ form a factorisable sequence of sheaves on $[X^n/S_n]$, then there exists a $S_n$-equivariant strict system of sheaves $\{F'_{\mathfrak A}\}_{\mathfrak A \vdash [n]}$ such that for every $\mathfrak A$ we have
	\[
	[F'_\mathfrak A] = [F_{\mathfrak A}]
	\]
	in $K([X^n/S_{\mathfrak A}])$.
\end{nprop}
\begin{proof}
	By Lemma \ref{thm:factorisableGivesSystem}, the sheaves $F_{\mathfrak A}$, where $\mathfrak A$ runs over partitions of $[n]$, underlie an $S_n$-equivariant system of sheaves on $X^n$.
	We write $(F,\phi)$ for this system.
	By Corollary \ref{thm:KTheoryMapSurjective}, there is an equality in $K(\Syst^S)$ of the form
	\[
	[(F,\phi)] = \sum_{j=1}^k i_j[\Psi(F'_j, \phi'_j)] - \sum_{j= k +1}^m i_j[\Psi(F'_j, \phi'_j)],
	\]
	where the $i_j > 0$ and $(F'_j,\phi'_j)$ are $S_n$-equivariant strict systems of sheaves.
	Taking
	\[
	F'_\mathfrak A = \left(\bigoplus_{j=1}^k (F'_j)^{\oplus i_j}_\mathfrak A \right)\oplus \left(\bigoplus_{j= k +1}^m (F'_j)^{\oplus i_j}_{\mathfrak A}\right)[1]
	\]
	gives the object we want.
\end{proof}

\section{Upgrading $G_n$ to a complex}
\label{sec:definingAComplex}
In this section, we let $(F_\mathfrak A, \phi_{\mathfrak A, \mathfrak B})$ be a strict $S_n$-equivariant system of sheaves on $X^n$, with respect to the stratification $(X^n_{\alpha})_{\alpha \in \mathcal I}$ explained in Example \ref{ex:mainExampleSystemContext}.

From the sheaves $F_{\mathfrak A}$ and their $S_n$-structure isomorphisms $\rho_{\mathfrak A, \sigma}$ we can define the sheaves $F_T, F_{S_nT}$ and $G_{n,k}$ exactly as in Section \ref{sec:explicitPlethysticLogarithm}.
The goal of this section is to show that the sheaf $G_n = \oplus_{k=0}^{n-1} G_{n,k}[k]$ can be equipped with a differential, built from the homomorphisms $\phi_{\mathfrak A, \mathfrak B}$, and to show that its cohomology vanishes away from the small diagonal of $X^n$ (Proposition \ref{thm:restrictionIsAcyclic}).

\begin{ndefn}
	Let $T$ be an index tree, and let $v$ be a non-leaf, non-root node.
	We define the \textbf{ordinary contraction} $O(T,v)$ to be the index tree where the node $v$ is deleted, and the parent node of $v$ is connected with all child nodes of $v$.
	The label sets of all the nodes in $O(T,v)$ are kept the same as in $T$.
\end{ndefn}

\begin{ndefn}
	In an index tree, we say that a node is \textbf{exceptional} if it is not a leaf and all its children are leaves.
\end{ndefn}

\begin{ndefn}
	Let $T$ be an index tree and $v$ an exceptional node of $T$. 
	We define the \textbf{exceptional contraction} $E(T,v)$ to be the index tree where all children of $v$ are deleted.
	The label sets of all nodes in $E(T,v)$ are kept the same as in $T$.
\end{ndefn}

We write $T \to T'$ to denote a contraction where $T' = O(T,v)$ or $T' = E(T,v)$ for some node $v$.
Note that for the set of non-leaf node labels $P(T)$, we have $P(T') = P(T) \setminus \{L_v\}$.
\begin{nremark}
	Identifying an index tree $T$ with its set of node labels $L(T) \subseteq 2^{[n]}$ as in Remark \ref{rmk:treesAsSets}, an ordinary contraction at $v$ is the operation of removing $L_v$ from $L(T)$.
	An exceptional contraction at $v$ is the operation of removing all $L_{v'}$ with $L_{v'} \subsetneq L_v$ from $L(T)$.
\end{nremark}

\begin{ndefn}
	Let $T$ be an index tree, and let $v$ be a non-leaf node of $T$.
	Define $s(v) = (-1)^l$, where $l$ is the number of $L \in P(T)$ such that $L < L_v$ in the binary ordering.
\end{ndefn}

\begin{ndefn}
	Given a non-leaf, non-root node $v$ of an index tree $T$, we define the ordinary contraction homomorphism $o(T,v) \colon F_{T} \to F_{O(T,v)}$ by
	\[
	o(T,v) \colon F_{T} = F_{\mathfrak A(T)} \overset{s(v)\id}\to F_{\mathfrak A(O(T,v))} = F_{O(T,v)}.
	\]
\end{ndefn}

\begin{ndefn}
	Given an exceptional node $v$ of an index tree $T$, we define the exceptional contraction homomorphism $e(T,v) \colon F_T \to F_{E(T,v)}$ by
	\[
	e(T,v) \colon F_{T} = F_{\mathfrak A(T)} \overset{-s(v)\phi_{\mathfrak A(T),\mathfrak A(E(T,v))}}{\to} F_{\mathfrak A(E(T,v))} = F_{E(T,v)}.
	\]
\end{ndefn}

Recall that
\[
G_{n,k} = \bigoplus_{T \in \mathcal T_{n,k}} F_{T}.
\]
We let
\[
d_{n,k} \colon G_{n,k} \to G_{n,k-1}
\]
be the sum over all $o(T,v)$ and $e(T,v)$ with $T \in \mathcal T_{n,k}$.
\begin{nlemma}
	The homomorphism $d_{n,k}$ is $S_n$-invariant.
\end{nlemma}
\begin{proof}
	Let $T$ be an index tree, let $v$ be a non-leaf, non-root node of $T$, and let $\sigma \in S_n$.
	Then in the commutative diagram
	\[
	\begin{tikzcd}[column sep = huge]
		(\sigma_X^{-1})^*F_T \arrow[r, "{(\sigma_X^{-1})^*o(T,v)}"] \arrow[d, "\rho_{T,\sigma}"] & (\sigma_X^{-1})^*F_{O(T,v)} \arrow[d, "\rho_{O(T,v),\sigma}"] \\
		F_{\sigma_X^{-1}(T)} \arrow[r, "{o(\sigma^{-1}T,\sigma^{-1}v)}"]  & F_{O(\sigma^{-1}T,\sigma^{-1}v)}
	\end{tikzcd}
\]
both compositions agree up to sign with the isomorphism 
\[
\rho_{\mathfrak A(T),\sigma} \colon (\sigma_X^{-1})^*F_{\mathfrak A(T)} \to F_{\sigma^{-1}\mathfrak A(T)}.
\]
Passing through the lower left hand corner first modifies this by a sign $(-1)^{l(T,\sigma)}$ from the definition of $\rho_{T,\sigma}$ in \eqref{eqn:rhoSigma}, and secondly by a sign $s(\sigma^{-1}(v))$ from the definition of $o(\sigma^{-1}T,\sigma^{-1}v)$.
The total sign is $(-1)^N$, where
\[
N = |\{A < B \in P(T) \mid \sigma^{-1}A > \sigma^{-1}B\}| + |\{A \in P(\sigma^{-1}T) \mid A < L_{\sigma^{-1}v}\}|
\]
Passing through the upper right hand corner introduces a sign $(-1)^{l({O(T,v)}, \sigma)}$, as well as a sign $(-1)^{s(v)}$, which amounts to the sign $(-1)^{N'}$, with
\[
N' = |\{A < B \in P(O(T,v)) \mid \sigma^{-1}A > \sigma^{-1}B\}| + |\{A \in P(T) \mid A < L_v\}|
\]
We get
\[
N' - N = 2|\{A \in P(T) \mid A > L_v \text{ and } \sigma^{-1}A < \sigma^{-1}L_v\}|,
\]
so $(-1)^N = (-1)^{N'}$ and the diagram commutes.

A similar argument shows that the $e(T,v)$ commute with the $\rho_{T,\sigma}$ in the same way, so $d_{n,k}$ is $S_n$-invariant.
\end{proof}

\begin{nlemma}
	We have $d_{n,k-1} \circ d_{n,k} = 0$.
\end{nlemma}
\begin{proof}
	Given two trees $T, T'$ with $k$ and $k-2$ non-leaf nodes, respectively, we must show that the sum of all contractions composing to a map $F_T \to F_{T'}$ are 0.
	We may assume that there are nodes $v_1, v_2 \in T$ which are not leaves, such that contracting at the $v_i$ in order will let us turn $T$ into $T'$.
	
	We have the following cases.
	\begin{enumerate}
		\item Neither $v_1$ nor $v_2$ is a descendant of each other.
		In this case we get
		\[
		o(O(T,v_1), v_2)o(T,v_1) + o(O(T,v_2), v_1)o(T, v_2) = 0
		\]
		and similarly for any exceptional contractions at $v_1$ or $v_2$.
		Thus the sum of all composed contraction homomorphisms from contracting at the $v_i$ is 0.
		\item The node $v_2$ is a descendant of $v_1$.
		There are the following subcases.
		\begin{enumerate}
			\item $v_2$ is not exceptional. 
			Then $v_1$ is not exceptional in either $T$ and $O(T, v_2)$, and $v_2$ is ordinary in $O(T,v_1)$, so the sum of all ways of contracting is
			\[
			o(O(T,v_1),v_2)o(T,v_1) + o(O(T,v_2), v_1)o(T,v_2) = 0
			\]
			as above.
			\item $v_2$ is exceptional, and $v_1$ has a non-leaf child not equal to $v_2$.
			Then $v_1$ is ordinary in $T, O(T,v_2)$ and $E(T,v_2)$, and so the sum of all ways of contracting is
			\begin{align*}
			&o(O(T,v_2),v_1)o(T,v_2)  + o(O(T,v_1),v_2)o(T,v_1) \\ 
			&+ e(O(T,v_1),v_2)o(T,v_1) + o(E(T,v_2),v_1)e(T,v_2).
			\end{align*}
			We see that the first two maps and the last two maps sum to 0.
			\item $v_2$ is exceptional, and $v_1$ has $v_2$ as its only non-leaf child.
			If $v_1$ is not the root, then we get 
			\begin{align*}
			&o(O(T,v_2),v_1)o(T,v_2) + o(O(T,v_1),v_2)o(T,v_1) \\
			&+ e(O(T,v_1),v_2)o(T,v_1) + o(E(T,v_2),v_1)e(T,v_2) \\
			& + e(O(T,v_2),v_1)o(T,v_2) + e(E(T,v_2), v_1)e(T,v_2).
			\end{align*}
			Here the terms on each line sum to 0.
			If $v_1$ is the root, then the first four terms do not appear.
		\end{enumerate}
	\end{enumerate}
\end{proof}

\subsection{Acyclicity of the complex}
\label{sec:acyclicity}
By the above lemmas, the homomorphism $d_n = \sum_{k=1}^{n-1} d_{n,k}$ defines a differential on the sheaf
\[
G_n = \bigoplus_{k=0}^{n-1} G_{n,k}[k]
\]
on $[X^n/S_n]$, i.e.~we have $d_n \colon G_n \to G_n[1]$ with $d_n^2 = 0$.
Let $\Delta \subseteq X^n$ be the small diagonal.
\begin{nprop}
	\label{thm:restrictionIsAcyclic}
	The cohomology sheaf $H^*(G_n) = \ker d_n/\im d_n$ is supported set-theoretically on $\Delta$.
\end{nprop}
We have
\[
X^n \setminus \Delta = \bigcup_{\mathfrak B} U_\mathfrak B.
\]
Here the union in $\mathfrak B$ runs over all 2-element set partitions $\mathfrak B = \{B_1, B_2\}$ of $[n]$, and $U_\mathfrak B$ is as in Section \ref{sec:systems}, concretely
\[
U_{\mathfrak B} = \{(x_i)_{i=1}^n \mid i \in B_1, j \in B_2 \Rightarrow x_i \not= x_j\}.
\]
Proposition \ref{thm:restrictionIsAcyclic} therefore follows from the following lemma.
\begin{nlemma}
	For every partition $\mathfrak B = \{B_1, B_2\}$ of $[n]$, we have
	\[
	H^*(G_n)|_{U_{\mathfrak B}} = 0.
	\]
\end{nlemma}
\begin{proof}
	In Lemma \ref{thm:fromPsiToFiltration}, we explain how a choice of a function $\psi \colon \cT_n \to S$ to an ordered set $S$ which decreases under contractions induces a filtration on the complex $G_n$.
	We then construct a specific such $\psi$, and show in Lemma \ref{thm:contractionsIncreasePsi} that it decreases under contractions.
	In Lemma \ref{thm:psiPreservingContractions} and Corollary \ref{thm:subquotientsAcyclic} we show that the associated subquotients of $G_n$ are acyclic over $U_{\mathfrak B}$, and it follows that $G_n$ is acyclic over $U_{\mathfrak B}$.
\end{proof}

\begin{nlemma}
	\label{thm:fromPsiToFiltration}
	Let $(S, \le)$ be a totally ordered set, and let $\psi \colon \mathcal T_{n} \to S$ be a map with the property that for every contraction $T \to T'$, we have $\psi(T) \ge \psi(T')$.
	
	Then there exists a filtration of the complex $G_n$ where every subquotient is of the form
	\[
	\bigoplus_{T \in \psi^{-1}(s)} F_T
	\]
	for some $s \in S$, with differential given by the sum of all differentials $F_T \to F_{T'}$ for contractions $T \to T'$ with $T, T' \in \psi^{-1}(s)$.
\end{nlemma}
\begin{proof}
	Since $\mathcal{T}_n$ is finite, we can map $\psi(\mathcal{T}_n) \subseteq S$ into $\ZZ$ in an order preserving way, and so we may assume that $S = \ZZ$ with the standard ordering.
	For each $i \in \ZZ$, define a subcomplex $G_{n}^i \subseteq G_n$ by
	\[
	\bigoplus_{T \in \psi^{-1}(-\infty, i]} F_T \subseteq \bigoplus_{T \in \mathcal T_n} F_T,
	\]
	where the condition that any contraction $T \to T'$ gives $\psi(T) \ge \psi(T')$ implies that $d_n(G^i_n) \subseteq G^i_n$.
	The subquotients $G^i_n/G^{i-1}_n$ then have the form claimed.
\end{proof}

Let $\le'$ be the total order on $2^{[n]}$ which is opposite to the binary total order, so that in particular $A \subseteq B \Rightarrow B \le' A$.
Let
\[
S = 2^{[n]} \times \NN^4,
\]
equipped with the lexicographical total ordering, where $2^{[n]}$ is ordered by $\le'$ and the $\NN$-factors in the usual way.
Define a map 
\[
\psi = (\psi_1, \psi_2, \psi_3, \psi_4, \psi_5) \colon \mathcal T_n \to S
\]
as follows.

Let $T$ be an index tree, and let 
\[
Q(T) = \{L_v \mid v\text{ a node of }T, L_v \not\subseteq B_1, L_v \not\subseteq B_2\}
\]
If $r$ is the root node of $T$, then $L_r = [n] \in Q(T)$, so $Q(T)$ is nonempty.

We let $d(T)$ be the node such that $L_{d(T)} \in Q(T)$ and this is the maximal element of $Q(T)$ with respect to $\le'$.
If $v \succ d(T)$, we then have
\[
v \succ d(T) \Rightarrow L_v \subsetneq L_{d(T)} \Rightarrow L_v >' L_{d(T)} \Rightarrow L_v \not\in Q(T),
\]
so that either $L_v \subseteq B_1 \cap L_{d(T)}$ or $L_v \subseteq B_2 \cap L_{d(T)}$.

We now define $\psi$ by
\begin{gather*}
	\psi_1(T) = L_{d(T)} \\
	\psi_2(T) = |\{v \text{ a node of } T \mid v\not\succeq d(T)\}| \\
	\psi_3(T) = |\{v \succ d(T) \mid L_v \subsetneq B_1 \cap L_{d(T)}\}| \\
	\psi_4(T) = \psi_3(T) \times |\{v \succ d(T) \mid L_v \subseteq B_2 \cap L_{d(T)}\}| \\
	\psi_5(T) = |\{v \succ d(T) \mid L_v \subsetneq B_2 \cap L_{d(T)}\}|
\end{gather*}

\begin{nlemma}
	\label{thm:contractionsIncreasePsi}
	If $T'$ is a contraction of $T$, then $\psi(T) \ge \psi(T')$.
\end{nlemma}
\begin{proof}
	We have $Q(T') \subseteq Q(T)$, which gives $\psi_1(T) \ge' \psi_1(T')$.
	Assume $\psi_1(T) = \psi_1(T')$, so that $L_{d(T)} = L_{d(T')}$.
	
	If the contraction is at a node $v \not\succeq d(T)$, then $\psi_2(T) > \psi_2(T')$.
	If the contraction is at $d(T)$, then the fact that $L_{d(T)} = L_{d(T')}$ implies the contraction must be exceptional, in which case
	we have $\psi_2(T) = \psi_2(T')$.
	If the contraction is at some node $v \succ d(T)$, we also get $\psi_2(T) = \psi_2(T')$.
	
	Assume now that $\psi_2(T) = \psi_2(T')$.
	If the contraction is exceptional at $d(T)$, we have
	\[
	\psi_3(T') = \psi_4(T') = \psi_5(T') = 0,
	\]
	so $\psi(T) \ge \psi(T')$.
	If the contraction happens at $v \succ d(T)$, then we get 
	\[
	\{L_v \mid v \succ d(T)\} \supsetneq \{L_{v'} \mid v' \succ d(T')\},
	\]
	from which we get
	\[
	\psi_{i}(T) \ge \psi_{i}(T')\ \ \ \ \ i = 3,4,5.
	\]
\end{proof}

Let us say that a contraction $T \to T'$ \textbf{preserves $\psi$} if $\psi(T) = \psi(T')$.
\begin{nlemma}
	\label{thm:psiPreservingContractions}
	For every index tree $T$ of order $n$, there is exactly one $\psi$-preserving contraction with $T$ as source or target.
	
	If the contraction $T \to T'$ preserves $\psi$, then the homomorphism $\phi_{T,T'} \colon F_T \to F_{T'}$ is an isomorphism over $U_\mathfrak B$.
\end{nlemma}
\begin{proof}
	If a contraction preserves $\psi$, then $\psi_1(T) = \psi_1(T')$ and $\psi_2(T) = \psi_2(T')$ imply that either \textit{(i)} the contraction is exceptional at $d(T)$, or \textit{(ii)} happens at a node $v \succ d(T)$.
	The equalities $\psi_3(T) = \psi_3(T')$ and $\psi_5(T) = \psi_5(T')$ imply that in case \textit{(i)}, there are exactly two descendants $v_1,v_2$ of $d(T)$, and these satisfy $L_{v_i} = B_i \cap L_{d(T)}$.
	In case \textit{(ii)}, the same equalities give that the contraction must be ordinary at a node $v$ with $L_v = B_i \cap L_{d(T)}$ for $i = 1$ or $2$.
	
	Let us now say a contraction $T \to T'$ is of type $P$ if it is exceptional at $d(T)$, and $d(T)$ has precisely two descendants $v_1,v_2$ with $L_{v_i} = B_i \cap L_{d(T)}$.
	We say it is of type $Q1$ (resp.~$Q2$) if it is ordinary at a node $v$ with $L_v = B_1 \cap L_{d(T)}$ (resp.~$L_v = B_2 \cap L_{d(T)}$).
	The previous paragraph shows that every contraction such that $\psi(T) = \psi(T')$ must be of type $P, Q1$ or $Q2$.
	
	For a tree $T$, we either have that $d(T)$ is a leaf, or else one of the following five subcases.
	\begin{enumerate}
			\item There is no node $v_1$ such that $L_{v_1} = B_1 \cap L_{d(T)}$.
			\item There is a node $v_1$ such that $L_{v_1} = B_1 \cap L_{d(T)}$, and:
			\begin{enumerate}
				\item $v_1$ is not a leaf.
				\item $v_1$ is a leaf, and:
				\begin{enumerate}
					\item There is no node $v_2$ such that $L_{v_2} = B_2 \cap L_{d(T)}$.
					\item There is a node $v_2$ such that $L_{v_2} = B_2 \cap L_{d(T)}$, and:
					\begin{enumerate}
						\item $v_2$ is not a leaf.
						\item $v_2$ is a leaf.
					\end{enumerate}
				\end{enumerate}				
			\end{enumerate}
	\end{enumerate}
We now analyse the possible contractions of type $P, Q1, Q2$ with source or target $T$ for each of these cases, and determine which satisfy $\psi(T) = \psi(T')$.
\begin{enumerate}
	\item[($d(T)$ is a leaf)] There is no contraction of type $P$, $Q1$ or $Q2$ from $T$.
	There are no contractions to $T$ of type $Q1$ or $Q2$.
	There is a unique contraction $T' \to T$ of type $P$, where $T'$ is obtained by adding two nodes $v_1, v_2$ to $T$ such that $L_{v_i} = B_i \cap L_{d(T)}$.
	This contraction preserves $\psi$.
	\item[(1)] There is no contraction from $T$ of type $P$ or $Q1$.
	There may be a contraction of type $Q2$.
	Since $L_{d(T)}$ is the union of $L_v$ over all $v \succ d(T)$, and for all these we have $L_v \subseteq B_i \cap L_{d(T)}$ for $i = 1$ or $2$, there must be at least one $v \succ d(T)$ with $L_v \subsetneq B_1 \cap L_{d(T)}$.
	Therefore $\psi_3(T) > 0$, and so $\psi_4$ decreases under a contraction of type $Q2$ and it is not $\psi$-preserving.
	
	There is no contraction to $T$ of type $P$.
	There is a unique contraction to $T$ of type $Q1$, and $\psi$ is preserved under this contraction.
	There may be a contraction to $T$ of type $Q2$, but as above $\psi_4$ decreases under this contraction.
	\item[(2a)] There is no contraction from $T$ of type $P$. There is a contraction from $T$ of type $Q1$, and possibly one of type $Q2$.
	The contraction of type $Q1$ preserves $\psi$, but the one of type $Q2$ does not, because $\psi_3(T) > 0$, and so $\psi_4$ becomes strictly smaller under this contraction.
	
	There is no contraction to $T$ of type $P$ or $Q1$.
	There may be a contraction to $T$ of type $Q2$, but as above $\psi_4$ decreases under this contraction, so it is not $\psi$-preserving.

\item[(2bi)] There is no contraction from $T$ of type $P$, $Q1$ or $Q2$.
There is no contraction to $T$ of type $P$ or $Q1$.
There is a unique contraction $T' \to T$ of type $Q2$.
Since we have a leaf $v_1$ with $L_{v_1} = B_1 \cap L_{d(T)}$, we have $\psi_3(T) = \psi_3(T') = 0$, which gives $\psi_4(T) = \psi_4(T') = 0$.
Clearly $\psi_5(T) = \psi_5(T')$, so this is $\psi$-preserving.
	
	\item[(2biiA)] There is no contraction from $T$ of type $P$ or $Q1$. 
	There is a contraction from $T$ of type $Q2$, and this preserves $\psi$.
	There is no contraction of type $P$, $Q1$ or $Q2$ to $T$.

	\item[(2biiB)] There is one contraction from $T$ of type $P$, which preserves $\psi$, and none of type $Q1,Q2$.
	There are no contractions to $T$ of type $P, Q1$ or $Q2$.
	
\end{enumerate}

This shows that every $T$ is the source or target of precisely one $\psi$-preserving contraction.

A contraction $T \to T'$ of type $Q1$ or $Q2$ gives $\mathfrak A(T) = \mathfrak A(T')$, so $\phi_{T,T'}$ is an isomorphism.
For a contraction $T \to T'$ of type $P$, we can write
\begin{gather*}
\mathfrak A(T) = \{A_1,\ldots, A_{l}, B_1 \cap L_{d(T)}, B_2 \cap L_{d(T)}\} \\
\mathfrak A(T') = \{A_1,\ldots, A_{l}, L_{d(T)}\}.
\end{gather*}
We see that $\mathfrak A(T) \sim_{\mathfrak B} \mathfrak A(T')$, and so $\phi_{T,T'}$ is an isomorphism over $U_{\mathfrak B}$.
\end{proof}

\begin{ncor}
	\label{thm:subquotientsAcyclic}
	Every subquotient of the filtration of $G_n$ associated with $\psi$ is acyclic on $U_{\mathfrak B}$.
\end{ncor}
\begin{proof}
	Each subquotient of $G_{n}$ from the filtration defined by $\psi$ is of the form
	\[
	\bigoplus_{T \in \psi^{-1}(x)} F_T[k_T],
	\]
	for some $x \in S$, where $k_T$ is the number of non-leaf nodes of $T$.
	By Lemma \ref{thm:psiPreservingContractions}, if $C$ is the set of contractions $T \to T'$ with $\psi(T) = \psi(T') = x$, this equals
	\[
	\bigoplus_{(T \to T') \in C} (F_T[k_T] \oplus F_{T'}[k_T-1]),
	\]
	and the differential is the sum of all $\phi_{T,T'}$.
	By Lemma \ref{thm:psiPreservingContractions} again, these are isomorphisms over $U_{\mathfrak B}$, and so the complex is acyclic.
\end{proof}

\section{Proof of Proposition \ref{thm:mainProposition}}
\label{sec:proof}
\begin{nprop}
	\label{thm:mainGeometricProposition}
	Let $(F_n)_{n=1}^\infty$ be a factorisable sequence of sheaves on $[X^n/S_n]$.
	Let $S_n$ act trivially on $X$, and let $i_n \colon [X/S_n] \to [X^n/S_n]$ be the inclusion of the small diagonal.
	There exists a sequence $(H_n)_{n=1}^\infty$ of sheaves on $[X/S_n]$, such that
	\[
	1 + \sum_{n=1}^\infty [F_n]q^n = \shExp\left(\sum_{n=1}^\infty(i_n)_*[H_n]q^n\right). 
	\]
\end{nprop}
\begin{proof}
The condition on the $H_n$ is equivalent to saying that for each $n$, we have $(i_n)_*[H_n] = [G_n]$, where
\[
G_n = \langle q^n \rangle \shLog\left( 1 + \sum_{n=1}^\infty F_n q^n \right).
\]
By Lemma \ref{thm:formulaForLog}, we have that 
\[
G_n = \bigoplus_{k=0}^{n-1} \left(\bigoplus_{S_nT \in \mathcal T_{n,k}/S_n} F_{S_nT}[k]\right),
\]
and by Remark \ref{remark:alternativeDefT}, the class $[G_n]$ is thus a function of the classes $[F_{\mathfrak A}]$ as $\mathfrak A$ runs over set partitions of $[n]$.
By Lemma \ref{thm:factorisableGivesSystem}, these sheaves $F_\mathfrak A$ gives rise to an $S_n$-equivariant system of sheaves on $X^n$ in the sense of Example \ref{ex:mainExampleSystemContext}.

By Proposition \ref{thm:systemToGoodSystem}, we may replace the system $(F_\mathfrak A)$ by a strict system $(F'_\mathfrak A)$ such that $[F_\mathfrak A] = [F'_{\mathfrak A}]$ for all $\mathfrak A$.
By this replacement we obtain a new sheaf $G'_n$ such that $[G'_n] = [G_n]$.

By Proposition \ref{thm:restrictionIsAcyclic}, we have that $G'_n$ is equipped with a differential such that $G'_n$ is acyclic after restriction to $X^n \setminus X$.
It follows by devissage that 
\[
[G_n] = [G'_n] = (i_n)_*[H_n]
\]
for some sheaf $H_n$ on $[X/S_n]$, which is what we wanted.
\end{proof}

\begin{proof}[Proof of Proposition \ref{thm:mainProposition}]
	Let $\pi \colon [X^n/S_n] \to \Sym^n X$ be the coarse moduli space map, and let $F_n = \pi^*\cF_n$, so the sequence $(F_n)_{n=1}^\infty$ of sheaves on $[X^n/S_n]$ is factorisable.
	
	Let $p$ denote the morphisms from $[X/S_n]$ and $[X^n/S_n]$ to $[\pt/S_n]$, and $Rp_*$ the corresponding pushforwards on $K_{\mathbb T}(-)_\loc$.
	The plethystic exponential commutes with this derived pushforward, in the sense that for any sequence $(\alpha_n)_{n=1}^\infty$ with $\alpha_i \in K_{\mathbb T}([X^i/S_i])$, we have
	\[
	Rp_*\left(\locExp\left(\sum_{n=1}^\infty \alpha_n q^n\right)\right) = \locExp\left(\left(\sum_{n=1}^\infty Rp_*(\alpha_n) q^n\right)\right)
	\]
	With $H_n$ the sheaves on $[X/S_n]$ obtained in Proposition \ref{thm:mainGeometricProposition}, we then get
	\begin{align*}
1 + \sum_{n=1}^\infty Rp_*[F_n]q^n &= Rp_*\left(\locExp\left(\sum_{n=1}^\infty (i_n)_*[H_n]q^n\right)\right) \\
	&= \locExp\left(\sum_{n=1}^\infty Rp_*((i_n)_*[H_n])q^n\right) \\
	&= \locExp\left(\sum_{n=1}^\infty Rp_*[H_n]q^n\right)
	\end{align*}
	Taking $S_n$-invariants and traces on both sides now yields, by Lemma \ref{thm:numericalLemma}, that
	\begin{align*}
	1 + \sum_{n=1}^\infty \chi(\Sym^n X, \cF_n)q^n =& 1 + \sum_{n=1}^\infty \chi(X^n, F_n)^{S_n} \\
	=& \tr_t\left(1 + \sum Rp_*[F_n]^{S_n}q^n\right) \\
	=& \Exp \left(\sum_{n=1}^\infty \tr_t(Rp_*[H_n]^{S_n})q^n\right)\\
	=& \Exp \left(\sum_{n=1}^\infty \chi(X,H_n^{S_n})q^n\right).
	\end{align*}
	Taking $\cE_n = H_n^{S_n}$ gives Proposition \ref{thm:mainProposition}.
\end{proof}
\printbibliography
\end{document}